\documentclass[a4paper, 10pt, reqno]{amsart}

\usepackage[T1]{fontenc}
\usepackage[utf8]{inputenc}
\usepackage{textcomp}
\usepackage{microtype}
\usepackage{lmodern}
\usepackage{amsmath}
\usepackage{amsthm}
\usepackage{amssymb}
\usepackage{amsfonts}
\usepackage{relsize}
\usepackage[mathscr]{eucal}
\usepackage[linktocpage]{hyperref}
\usepackage{aliascnt}
\usepackage[sort,capitalize]{cleveref}
\usepackage{enumitem}			
\usepackage[dvipsnames]{xcolor}
\usepackage[msc-links, abbrev, non-sorted-cites]{amsrefs}
\usepackage{todonotes}
\usepackage{caption}

\SetMathAlphabet{\mathsf}{normal}{OT1}{lmss}{m}{n}
\SetMathAlphabet{\mathsf}{bold}{OT1}{lmss}{bx}{n}

\usepackage{tikz}

\makeatletter
\def\nonumberfootnote{\xdef\@thefnmark{}\@footnotetext}			
\makeatother

\definecolor{colorred}{HTML}{B00000}
\definecolor{colorgreen}{HTML}{258300}
\definecolor{colorblue}{HTML}{2e32fa}
\definecolor{coloryellow}{HTML}{cbbb1a}
\hypersetup{colorlinks=true, linkcolor=colorred, citecolor=colorgreen, urlcolor=coloryellow, pdfusetitle=true}

\numberwithin{equation}{section}

\newcommand{\rmc}{{\ensuremath{\mathrm{c}}}}
\newcommand{\rmd}{{\ensuremath{\mathrm{d}}}}

\newcommand{\sfa}{{\ensuremath{\mathsf{a}}}}
\newcommand{\sfb}{{\ensuremath{\mathsf{b}}}}

\newcommand{\scrC}{{\ensuremath{\mathscr{C}}}}
\newcommand{\scrD}{{\ensuremath{\mathscr{D}}}}

\newcommand{\scrK}{{\ensuremath{\mathscr{K}}}}

\newcommand{\N}{\boldsymbol{\mathrm{N}}}						
\newcommand{\R}{\boldsymbol{\mathrm{R}}}						
\renewcommand{\S}{\boldsymbol{\mathrm{S}}}						
\renewcommand{\d}{\,\mathrm{d}}				

\DeclareMathOperator*{\argmin}{argmin}		
\DeclareMathOperator{\supp}{spt}			

\theoremstyle{definition}

\newcommand{\newbumptheorem}[2]{%
  \newaliascnt{#1}{bump}%
  \newtheorem{#1}[#1]{#2}%
  \aliascntresetthe{#1}}

\theoremstyle{plain}
\newbumptheorem{theorem}{Theorem}
\newbumptheorem{proposition}{Proposition}
\newbumptheorem{definition}{Definition}
\newbumptheorem{lemma}{Lemma}
\newbumptheorem{corollary}{Corollary}
\newbumptheorem{hypothesis}{Hypothesis}
\newbumptheorem{idea}{Idea}
\newbumptheorem{conjecture}{Conjecture}
\newbumptheorem{question}{Question}

\theoremstyle{remark}
\newbumptheorem{remarkx}{Remark}
\newbumptheorem{examplex}{Example}
\newenvironment{remark}[1][]%
  {\pushQED{\qed}%
   \if\relax\detokenize{#1}\relax\remarkx\else\remarkx[#1]\fi}%
  {\popQED\endremarkx}
  {\pushQED{\qed}%
   \if\relax\detokenize{#1}\relax\examplex\else\examplex[#1]\fi}%
  {\popQED\endexamplex}
\crefname{remarkx}{Remark}{Remarks}
\crefname{examplex}{Example}{Examples}
\newbumptheorem{notation}{Notation}
\newbumptheorem{exercise}{Exercise}

\crefname{theorem}{Theorem}{Theorems}
\crefname{proposition}{Proposition}{Propositions}
\crefname{definition}{Definition}{Definitions}
\crefname{lemma}{Lemma}{Lemmas}
\crefname{corollary}{Corollary}{Corollaries}
\crefname{hypothesis}{Hypothesis}{Hypotheses}
\crefname{remark}{Remark}{Remarks}
\crefname{example}{Example}{Examples}
\crefname{notation}{Notation}{Notations}
\crefname{idea}{Idea}{Ideas}
\crefname{conjecture}{Conjecture}{Conjectures}

\crefformat{section}{{§}#2#1#3}
\crefformat{subsection}{{§}#2#1#3}
\crefformat{subsubsection}{{§}#2#1#3}
\crefformat{appendix}{{§}#2#1#3}

\crefmultiformat{theorem}{Theorems #2#1#3}{ and #2#1#3}{, #2#1#3}{, and #2#1#3}
\crefmultiformat{proposition}{Propositions #2#1#3}{ and #2#1#3}{, #2#1#3}{, and #2#1#3}
\crefmultiformat{definition}{Definitions #2#1#3}{ and #2#1#3}{, #2#1#3}{, and #2#1#3}
\crefmultiformat{exercise}{Exercises #2#1#3}{ and #2#1#3}{, #2#1#3}{, and #2#1#3}
\crefmultiformat{lemma}{Lemmas #2#1#3}{ and #2#1#3}{, #2#1#3}{, and #2#1#3}
\crefmultiformat{corollary}{Corollaries #2#1#3}{ and #2#1#3}{, #2#1#3}{, and #2#1#3}
\crefmultiformat{hypothesis}{Hypotheses #2#1#3}{ and #2#1#3}{, #2#1#3}{, and #2#1#3}
\crefmultiformat{remark}{Remarks #2#1#3}{ and #2#1#3}{, #2#1#3}{, and #2#1#3}
\crefmultiformat{example}{Examples #2#1#3}{ and #2#1#3}{, #2#1#3}{, and #2#1#3}
\crefmultiformat{notation}{Notations #2#1#3}{ and #2#1#3}{, #2#1#3}{, and #2#1#3}
\crefmultiformat{conjecture}{Conjectures #2#1#3}{ and #2#1#3}{, #2#1#3}{, and #2#1#3}

\crefmultiformat{section}{{§§}#2#1#3}{ and #2#1#3}{, #2#1#3}{, and #2#1#3}
\crefmultiformat{subsection}{{§§}#2#1#3}{ and #2#1#3}{, #2#1#3}{, and #2#1#3}
\crefmultiformat{subsubsection}{{§§}#2#1#3}{ and #2#1#3}{, #2#1#3}{, and #2#1#3}

\crefrangeformat{equation}{#3\textcolor{black}{(}#1\textcolor{black}{)}#4 to #5\textcolor{black}{(}#2\textcolor{black}{)}#6}

\newcommand{\Leb}{\mathscr{L}}				
\newcommand{\vol}{\mathrm{vol}}				

\newcommand{\Id}{\mathrm{Id}}				

\newcommand{\comp}{\rmc}						

\newcommand{\loc}{\mathrm{loc}}				

\DeclareMathOperator{\diam}{diam}			

\newcommand{\Len}{\mathrm{Len}}

\usepackage{blindtext}

\allowdisplaybreaks

\let\oldtocsection=\tocsection

\let\oldtocsubsection=\tocsubsection

\let\oldtocsubsubsection=\tocsubsubsection

\renewcommand{\tocsection}[2]{\hspace{0em}\oldtocsection{#1}{#2}}
\renewcommand{\tocsubsection}[2]{\hspace{1em}\oldtocsubsection{#1}{#2}}
\renewcommand{\tocsubsubsection}[2]{\hspace{2em}\oldtocsubsubsection{#1}{#2}}

\newcommand{\nocontentsline}[3]{}
\newcommand{\tocless}[2]{\bgroup\let\addcontentsline=\nocontentsline#1{#2}\egroup}

\newcommand{\mres}{\mathbin{\vrule height 1.6ex depth 0pt width 0.13ex\vrule height 0.13ex depth 0pt width 1.3ex}}

\allowdisplaybreaks

\makeatletter
\@namedef{subjclassname@2020}{\textup{2020} Mathematics Subject Classification}
\makeatother

\providecommand{\scrC}{\mathscr{C}}\providecommand{\Id}{\mathrm{Id}}

\begin{document}

\title[$p$-Harmonic gauge and spacetime isometries]{$p$-Harmonic gauge and the regularity of spacetime isometries}
\author{Mathias Braun}
\address{Institute of Mathematics, EPFL, 1015 Lausanne, Switzerland}
\email{\href{mailto:mathias.braun@epfl.ch}{mathias.braun@epfl.ch}}

\thanks{I am supported by the EPFL through a Bernoulli Instructorship. I sincerely thank Georgios Moschidis for enlightening discussions.}

\begin{abstract} We show every Lorentzian metric of class $C^{k,\alpha}$, where $k\in\N_0$ and $\alpha\in(0,1)$, admits local coordinates of class $C^{k+1,\alpha}$ all of whose coordinate functions are temporal functions that are harmonic with respect to a nonlinear $p$-wave operator, where $p\in(-\infty,1)$. Unlike the classical wave operator, this operator is elliptic on temporal functions. As an application, we show that an isometry between Lorentzian metrics of class $C^{k,\alpha}$ is of class $C^{k+1,\alpha}$.
\end{abstract}

\maketitle

\thispagestyle{empty}

\tableofcontents

\addtocontents{toc}{\protect\setcounter{tocdepth}{2}}

\section{Introduction}

\subsection{Harmonic gauge} In Riemannian geometry, harmonic coordinates yield optimal regularity, as pointed out by DeTurck--Kazdan \cite{deturck-kazdan1981}. They proved that if a Riemannian metric is of class $C^{k,\alpha}$ in some coordinates, where $k\in\N$ and $\alpha\in(0,1)$, then it is of class $C^{k,\alpha}$ in harmonic coordinates and the coordinate functions themselves are of class $C^{k+1,\alpha}$. This gain of one derivative, a consequence of elliptic regularity, is the mechanism behind several classical results at low regularity. One of them concerns isometries, cf.~\cref{Sub:ISometries}. Harmonic coordinates of class $C^{1,\alpha}$ still exist if the metric is merely of class $C^{0,\alpha}$, cf.~e.g.~Julin--Liimatainen--Salo \cite{julin-liimatainen-salo2017}. Here, Hölder continuity cannot be replaced by continuity: by Hartman--Wintner \cite{hartman-wintner1953}, there are continuous two-dimensional Riemannian metrics which admit no isothermal coordinates of class $C^1$ and for which the Dirichlet problem with nonconstant boundary data possesses no $C^1$-solution. In particular, such metrics admit no harmonic coordinates of regularity $C^1$.

In Lorentzian signature, harmonic coordinates (or wave coordinates) appear prominently as well. First used by Einstein \cite{einstein1916-integration} in a special case and subsequently, for the full field equations, by de Donder \cite{dedonder1921} and Lanczos \cite{lanczos1922}, they underlie the local well-posedness theory of the Einstein equations since Choquet-Bruhat \cite{choquet-bruhat1952}. Wave coordinates solve a hyperbolic PDE. Thus, they are determined by initial rather than boundary data and their regularity is measured by energy estimates, cf.~e.g.~Klainerman--Rodnianski \cite{klainerman-rodnianski2005}; in the Hölder scale, no derivative is gained. For merely Hölder Lorentzian metrics, the Cauchy problem for the wave equation need not be solvable even for smooth initial data, cf.~Colombini--De Giorgi--Spagnolo \cite{colombini-degiorgi-spagnolo1979}, and the standard construction of wave coordinates breaks down.

\subsection{$p$-Harmonic gauge} For Riemannian metrics, a nonlinear variant of harmonic coordinates was introduced by Liimatainen--Salo \cite{liimatainen-salo2014,liimatainen-salo2016} with Julin \cite{julin-liimatainen-salo2017}. They constructed $p$-harmonic coordinates, in \cite{julin-liimatainen-salo2017} for every Hölder Riemannian metric, where $p\in (1,\infty)$, and used them, for $p$ equal to the dimension, to prove regularity of conformal maps and local conformal flatness.

In this note, we introduce coordinates for Lorentzian metrics (possibly of low regularity) that solve an elliptic PDE, despite the hyperbolic nature of Lorentzian geometry. We call them \textbf{$p$-harmonic coordinates}: we replace the classical wave operator by the nonlinear $p$-wave operator
\begin{align*}
    \Box_pu := \vert\!\det g\vert^{-1/2}\,\partial_\alpha\big[\vert\!\det g\vert^{1/2}\,\big\vert\rmd u\big\vert_g^{p-2}\,g^{\alpha\beta}\,\partial_\beta u\big]
\end{align*}
on temporal functions (functions with future-directed timelike differential), where $p\in (-\infty,1)$ and
\begin{align*}
    \vert \rmd u\vert_g := \big[{-g^{-1}(\rmd u,\rmd u)}\big]^{1/2}.
\end{align*}
This operator was introduced by Mondino--Suhr \cite{mondino-suhr2022}. Its ellipticity on temporal functions was observed by Beran et al.~\cite{beran-braun-calisti-gigli-mccann-ohanyan-rott-samann+-}; see \cref{Re:Whyelliptic} for an illustration that also explains the somewhat unusual range of $p$. Braun et al.~\cite{braun-gigli-mccann-ohanyan-samann+,braun-gigli-mccann-ohanyan-samann++} used it to give elliptic proofs and extensions of the Lorentzian splitting theorems of Eschenburg \cite{eschenburg1988}, Galloway \cite{galloway1989-splitting}, and Newman \cite{newman1990}. Given ellipticity, the regularity of $p$-harmonic functions is governed by the classical theory of Schauder and De Giorgi--Nash--Moser. We refer to \cref{Sec:Prelim} for precise definitions about the nonlinear $p$-wave operator and $p$-harmonic temporal functions.

Throughout, let $d\in\N$ and $M\subset \R^{1+d}$ be open. Unless stated otherwise, local coordinate expressions, and in particular the regularity of a metric, are understood with respect to the identity chart; since our main results are local, $M$ should be regarded as a chart domain of a spacetime in applications. All Lorentzian metrics are assumed to be time-oriented. The range of Greek indices is $\{0,\dots,d\}$. Repeated indices, one as a subscript and one as a superscript in the same term, are summed over. Our signature convention for Lorentzian metrics is $-,+,\dots,+$.

Our first result shows $p$-harmonic temporal functions gain one derivative over the metric in the Hölder scale, as harmonic functions do in Riemannian geometry. Notably, this includes the case of merely Hölder metrics. 

\begin{proposition}[Regularity of $p$-harmonic temporal functions]\label{thm:reg} Let $g$ be a given Lorentzian metric on $M$ of class $\smash{C^{k,\alpha}_\loc}$, where $\smash{k\in\N_0}$ and $\alpha\in(0,1)$. Then every $p$-harmonic temporal function on $M$ is of class $\smash{C^{k+1,\alpha}_\loc}$.
\end{proposition}

When $g$ is smooth, using classical De Giorgi--Nash--Moser theory, Ohanyan \cite{ohanyan2026+} has recently shown smoothness of $p$-harmonic temporal functions. 

If $k\in\N$, the conclusion of \cref{thm:reg} follows from the classical method of difference quotients and Schauder estimates. For $k=0$, difference quotients of the coefficients are unbounded; we use a perturbation argument instead. For locally Lipschitz metrics, which corresponds to the endpoint $\alpha=1$ not covered by \cref{thm:reg} for $k=0$, our arguments entail $p$-harmonic temporal functions have second derivatives in $\smash{L^2_\loc}$, cf.~\cref{lem:H2}.

Since the cone of future-directed timelike covectors is open, it contains bases of the cotangent space. This allows us to use temporal functions for \emph{all} coordinate directions, which is our second result.

\begin{theorem}[$p$-Harmonic coordinates]\label{thm:coord} Let $g$ be a Lorentzian metric on $M$ of class $\smash{C_\loc^{k,\alpha}}$, where $k\in\N_0$ and $\alpha\in(0,1)$. Let $x_0\in M$. Then there exist $r>0$ and $p$-harmonic temporal functions $\smash{\tau^0,\dots,\tau^d\colon B_r(x_0)\to\R}$ such that $\smash{\Phi := (\tau^0,\dots,\tau^d)}$ is a $\smash{C^{k+1,\alpha}}$-diffeomorphism from $B_r(x_0)$ onto its image; in particular, the coefficients of $\smash{\Phi_*g}$ are of class $\smash{C^{k,\alpha}}$.

In particular, if $g$ is smooth, it admits smooth $p$-harmonic coordinates around each point in $M$.
\end{theorem}

In the coordinates $\smash{\tau^0,\dots,\tau^d}$ from \cref{thm:coord}, we have $\smash{\vert\rmd \tau^\mu\vert_g^2 = -g^{\mu\mu}}$, where $\smash{g^{\alpha\beta}}$ denote the coefficients of the inverse metric with respect to these coordinates. Hence, the metric satisfies $g^{\mu\mu}<0$ and the distributional gauge condition
\begin{align}\label{Eq:Gauge}
    \partial_\alpha\big[\vert\!\det g\vert^{1/2}\,\big[{-g^{\mu\mu}}\big]^{(p-2)/2}\,g^{\alpha\mu}\big] = 0
\end{align}
without summation over $\mu$. It replaces the classical harmonic gauge condition
\begin{align*}
    \partial_\alpha\big[\vert\!\det g\vert^{1/2}\,g^{\alpha\mu}\big]=0
\end{align*}
formally corresponding to $p=2$, which is excluded from our range of $p$.

We point out that \cref{thm:coord} contains a \emph{local} existence result for $p$-harmonic temporal functions. We obtain it by truncating the equation outside a compact subset of the timelike cone, solving the resulting Dirichlet problem for a strongly monotone operator, and verifying a posteriori that the differential of the solution stays in the untruncated region, cf.~\cref{prop:exist}. The truncation is required since the nonlinear $p$-wave operator degenerates as the differential approaches the light cone and ceases to be elliptic beyond it, cf.~\cref{Re:Whyelliptic}; by contrast, the Riemannian $p$-Laplacian degenerates only where the gradient vanishes. To the best of our knowledge, \cref{prop:exist} is the first local existence result for $p$-harmonic temporal functions which requires neither curvature nor symmetry assumptions. In Minkowski spacetime, $p$-harmonic temporal functions arise as minimizers of the Lorentzian $p$-energy on suitable domains, which are shown to be temporal and locally Lipschitz by De Philippis et al.~\cite{dephilippis-gigli-mccann-ohanyan-picerni-zanardini+}. \emph{Global} existence is a different matter. In the elliptic proofs of the Lorentzian splitting theorems by Braun et al.~\cite{braun-gigli-mccann-ohanyan-samann+,braun-gigli-mccann-ohanyan-samann++}, a globally $p$-harmonic function is provided by the Busemann function of a timelike line, whose $p$-harmonicity relies on timelike Ricci nonnegativity. Without a line, this construction is not available, cf.~Galloway--McCann--Ohanyan \cite{galloway-mccann-ohanyan2026+}. In view of \cite{braun-gigli-mccann-ohanyan-samann+,braun-gigli-mccann-ohanyan-samann++}, it is then natural to search for globally $p$-harmonic temporal functions directly; their existence is equivalent to Bartnik's splitting conjecture \cite{bartnik1988} on the spacetimes it concerns, see Ohanyan \cite{ohanyan2026+}.

\subsection{Regularity of spacetime isometries}\label{Sub:ISometries} As an application of our previous results, we revisit the question of regularity of spacetime isometries.

By a classical theorem of Myers--Steenrod \cite{myers-steenrod1939} and Palais \cite{palais1957}, distance-preserving maps between smooth Riemannian manifolds are smooth. Their proofs, as well as several results for nonsmooth metrics, rest on geodesics. Hartman \cite{hartman1952} showed a $C^1$-map that transforms a two-dimensional Riemannian metric of class $C^1$ into another one is of class $C^2$. At this regularity, the geodesic equation is already delicate: Hartman--Wintner \cite{hartman-wintner1951} pointed out its solutions are neither unique nor necessarily minimizing, even if the Christoffel symbols are Hölder continuous. In \cite{hartman1958}, Hartman removed the restrictions on dimension and signature by a different method, which avoids geodesics and treats the isometry condition as an overdetermined first-order system. Calabi--Hartman \cite{calabi-hartman1970} observed that the argument of \cite{myers-steenrod1939} is not correct for metrics of class $C^1$, since it employs normal coordinates, which need not exist at this regularity, and gave a new proof. For Hölder metrics, where no geodesic equation is available, they based their study of isometries on a claim about the regularity of minimizing curves. This claim, in turn, was falsified by Lytchak--Yaman \cite{lytchak-yaman2006}. The statement on isometries is true nevertheless and alternative proofs have been found: that isometries between Riemannian metrics of class $C^{k,\alpha}$ are of class $C^{k+1,\alpha}$ is due to Reshetnyak \cite{reshetnyak1978} for $k=0$, who also pointed out the gap in \cite{calabi-hartman1970}, to Shefel' \cite{shefel1982} for $k\geq 1$, and to Sabitov \cite{sabitov1993} in borderline cases. Taylor \cite{taylor2006} gave a unified proof by means of harmonic coordinates, which are preserved by isometries; in dimension two, for analytic surfaces and isothermal coordinates, this argument goes back to Wintner \cite{wintner1952}. We also refer to Matveev--Troyanov \cite{matveev-troyanov2017} for extensions to Finsler Lagrangians.

The preceding results do not extend to merely continuous Riemannian metrics. As pointed out by Calabi--Hartman \cite{calabi-hartman1970}, distance-preserving maps between such metrics need not be of regularity $C^1$. 

The $p$-wave operator is invariant under isometries. Together with \cref{thm:reg,thm:coord}, this yields our third result.

\begin{theorem}[Regularity of $C^1$-isometries]\label{thm:iso} Let $g$ and $g'$ be Lorentzian metrics of class $\smash{C^{k,\alpha}_\loc}$ on open subsets $\smash{M,M'\subset\R^{1+d}}$, respectively, where $k\in\N_0$ and $\alpha\in(0,1]$. Let $\psi\colon M\to M'$ be a $C^1$-isometry between $g$ and $g'$. Then $\psi$ is of class $\smash{C^{k+1,\alpha}_\loc}$.
\end{theorem}

For $k\in\N$, \cref{thm:iso} follows from a theorem of Hartman \cite[Thm.~II]{hartman1958}: every $C^1$-map which transforms a nonsingular symmetric matrix field of class $C^k$ into another one is of class $C^{k+1}$. For metrics of class $\smash{C^{k,\alpha}}$, $k\in\N$ and $\alpha\in(0,1]$, the transformation rule of the Christoffel symbols, which expresses the second derivatives of $\psi$ through $g$, $g'$, their first derivatives, and the first derivatives of $\psi$, then shows $\psi$ is of class $\smash{C^{k+1,\alpha}}$, cf.~Calabi--Hartman \cite[§5]{calabi-hartman1970}. For $k=0$ and $\alpha\in(0,1)$, the Christoffel symbols are not defined. This case is our contribution: Hartman's theorem is trivial here, as it asserts that a map which is assumed to be of class $C^1$ is of class $C^1$, whereas \cref{thm:iso} establishes Hölder continuity of its first derivatives. It is the Lorentzian counterpart of the theorem of Reshetnyak \cite{reshetnyak1978}; we are not aware of a previous proof. For $\alpha\in(0,1)$, the proof we give is uniform in $k$. Moreover, our approach via $p$-harmonic coordinates is fundamentally different from that of \cite{calabi-hartman1970}, which is based on distances and geodesics. The endpoint case $k=0$ and $\alpha=1$ of Lipschitz metrics, where \cref{thm:iso} shows isometries are of class $\smash{C^{1,1}}$, is new for Lorentzian metrics as well, to the best of our knowledge. Here the Christoffel symbols are defined only a.e.~and, once we have used $p$-harmonic coordinates to show the second derivatives of $\psi$ are in $\smash{L_\loc^2}$, the last step of our proof adapts the argument of Calabi--Hartman \cite[§5]{calabi-hartman1970} to this setting.

The hypothesis that $\psi$ is of class $C^1$ in \cref{thm:iso} can be relaxed, provided $\psi$ is compatible with the time orientations.

\begin{theorem}[Regularity of bi-Lipschitz isometries]\label{thm:bilip} Let $g$ and $g'$ be Lorentzian metrics of class $\smash{C^{k,\alpha}_\loc}$ on open subsets $\smash{M,M'\subset\R^{1+d}}$, respectively, where $k\in\N_0$ and $\alpha\in(0,1]$. Let $\psi\colon M\to M'$ be a locally bi-Lipschitz map that is a time-orientation-preserving isometry between $g$ and $g'$, cf.~\cref{Def:Isometry}. Then $\psi$ is of class $\smash{C^{k+1,\alpha}_\loc}$; in particular, $\psi^*g'=g$ everywhere on $M$.
\end{theorem}

Without local injectivity, the conclusion fails, as shown on Minkowski spacetime by the map $\smash{\psi\colon \R^{1+d}\to\R^{1+d}}$ given by $\smash{\psi(t,x^1,\dots,x^d) := (t,\vert x^1\vert,x^2,\dots,x^d)}$, which obeys all other hypotheses of \cref{thm:bilip}. Neither can the Hölder continuity of the metrics be dropped in general. Let $\smash{f\colon U\to\R^2}$ denote the map of Calabi--Hartman \cite[§6]{calabi-hartman1970}, defined on a neighborhood $\smash{U\subset\R^2}$ of the origin; it is a homeomorphism onto its image and an isometry from the Euclidean metric to a continuous Riemannian metric $h$. It is not of class $C^1$ but easily seen to be locally bi-Lipschitz. The map $\smash{\psi\colon\R\times U\to\R\times f(U)}$ given by $\smash{\psi(t,x) := (t,f(x))}$ pulls the continuous Lorentzian metric $-\rmd t^2+h$ back to the Minkowski metric outside the $t$-axis and satisfies all other hypotheses of \cref{thm:bilip}. Yet, $\psi$ fails to be differentiable along the $t$-axis.

The two a.e.~conditions in \cref{thm:bilip} are implied by the preservation of the time separation function; see \cref{Sub:RegSptIso} for the relevant notions.

\begin{theorem}[Regularity of maps preserving the time separation]\label{thm:tau} Let $g$ and $g'$ be strongly causal Lorentzian metrics of class $\smash{C^{k,\alpha}_\loc}$ on open subsets $\smash{M,M'\subset\R^{1+d}}$, respectively, where $k\in\N_0$ and $\alpha\in(0,1]$. Moreover, let $\psi\colon M\to M'$ be a locally bi-Lipschitz homeomorphism that preserves the respective time separation functions, cf.~\cref{Def:TimeSepPres}. Then $\psi$ is of class $\smash{C^{k+1,\alpha}_\loc}$ and satisfies $\psi^*g'=g$ on $M$.
\end{theorem}

We do not know whether the assumption that $\psi$ is locally bi-Lipschitz in \cref{thm:tau} can be dropped; for smooth metrics, this is the case by the theorem of Hawking--King--McCarthy \cite{hawking-king-mccarthy1976}.

\subsection{Remarks}
\begin{enumerate}[label=\textnormal{\arabic*.}]
    \item Below Lipschitz regularity, the causality theory of Lorentzian metrics degenerates; for instance, light cones may admit bubbles, cf.~Chruściel--Grant \cite{chrusciel-grant2012}. Apart from \cref{thm:tau}, our arguments are local and elliptic and are not affected by this.
    \item The Riemannian analog of \cref{thm:coord} for $p$ equal to the dimension or in $(1,\infty)$, respectively, is due to Liimatainen--Salo \cite{liimatainen-salo2014,liimatainen-salo2016} with Julin \cite{julin-liimatainen-salo2017}. In their application to conformal flatness, the conformally invariant exponent equal to the dimension plays a distinguished role. In Lorentzian signature, this exponent is $1+d$ and lies outside the elliptic range $(-\infty,1)$.
    \item The gauge condition \eqref{Eq:Gauge} states that the coordinate functions solve an elliptic equation. We will discuss applications to the Einstein equations in forthcoming work \cite{braun+}.
\end{enumerate}

\section{Nonlinear $p$-wave operator and $p$-harmonic temporal functions}\label{Sec:Prelim}

Throughout this section, $g$ is a Lorentzian metric on $M$ with continuous coefficients. We assume it is time-oriented by a continuous vector field $X$ on $M$ with $g(X,X)<0$ on $M$. The Euclidean norm of vectors, covectors, and matrices is denoted by $\vert\cdot\vert$. As usual, $B_r(x)$ is the open Euclidean ball of radius $r>0$ around $\smash{x\in\R^{1+d}}$. Let $\vol$ denote the volume measure induced by $g$. Fix $p\in(-\infty,1)$; all constants are allowed to depend on $p$ and $d$ without further notice.

Let $x\in M$. We call a nonzero covector $\smash{\zeta\in T_x^*M}$ \textbf{timelike} if $\smash{g^{-1}(x)(\zeta,\zeta)<0}$, \textbf{causal} if $\smash{g^{-1}(x)(\zeta,\zeta)\leq 0}$, and \textbf{future-directed} if $\zeta(X(x))>0$. A nonzero tangent vector $v\in T_xM$ is called timelike if $g(x)(v,v)<0$, causal if $g(x)(v,v)\leq 0$, and future-directed if $g(x)(v,X(x))<0$. With a slight abuse of notation, we will write $\smash{\vert \cdot\vert_g}$ both for $\smash{[-g(\cdot,\cdot)]^{1/2}}$ and $\smash{[-g^{-1}(\cdot,\cdot)]^{1/2}}$. In the sequel, we sometimes identify all cotangent spaces of $M$ with $\smash{\R^{1+d}}$ by means of the identity chart; accordingly, an element of $\smash{\R^{1+d}}$ can be regarded as a covector at each point of $M$.

\begin{definition}[Temporal function]\label{Def:Temporal} A function $u\colon M\to\R$ is called \textbf{temporal} if it is locally Lipschitz and for every compact $K\subset M$, there exists $c>0$ with
\begin{align}\label{Eq:LowBo}
    \min\{-g^{-1}(\rmd u,\rmd u), \rmd u(X)\}\geq c\quad\vol\mres K\textnormal{\textit{-a.e.}}
\end{align}
\end{definition}

\begin{definition}[Nonlinear $p$-wave operator] The \textbf{nonlinear $p$-wave operator} $\Box_p$ assigns to every temporal function $u\colon M\to\R$ the distribution
\begin{align*}
    \Box_pu := \vert\!\det g\vert^{-1/2}\,\partial_\alpha\big[\vert\!\det g\vert^{1/2}\,\big\vert\rmd u\big\vert_g^{p-2}\,g^{\alpha\beta}\,\partial_\beta u\big],
\end{align*}
understood with respect to $\smash{\vol}$, i.e.~for every $\smash{\varphi\in C_\comp^\infty(M)}$,
\begin{align}\label{Eq:Weak}
    \big\langle\Box_pu,\varphi\big\rangle := -\int_M\big\vert\rmd u\big\vert_g^{p-2}\,g^{-1}(\rmd u,\rmd\varphi)\d\vol.
\end{align}
\end{definition}

It is clear that this definition makes sense for merely continuous $g$.

This nonlinear $p$-wave operator has first appeared in Mondino--Suhr \cite{mondino-suhr2022}. Its ellipticity has first been observed by Beran et al.~\cite{beran-braun-calisti-gigli-mccann-ohanyan-rott-samann+-}; this explains our choice of $p$, illustrated in \cref{Re:Whyelliptic}. It was then employed to simplify and generalize the Lorentzian splitting theorems of Eschenburg \cite{eschenburg1988}, Galloway \cite{galloway1989-splitting}, and Newman \cite{newman1990} by Braun et al.~\cite{braun-gigli-mccann-ohanyan-samann+,braun-gigli-mccann-ohanyan-samann++}; see also the Finslerian generalization by Caponio--Ohanyan--Ohta \cite{caponio-ohanyan-ohta2024+}. Exact representation formulas for this operator were shown by Braun \cite{braun2024+}. We refer to the surveys of McCann \cite{mccann+} and Braun \cite{braun2025}. The nonlinear $p$-wave operator is properly associated to a $p$-Hamiltonian introduced independently by McCann \cite{mccann2020}, Minguzzi \cite{minguzzi2019-causality}, and Mondino--Suhr \cite{mondino-suhr2022} that is strictly convex on the future time cone (which already hints at ellipticity, cf.~\cref{lem:ell}) in the sense of Hamiltonian geometry, cf.~Ohta \cite{ohta2014}.

\begin{remark}[Admissible test functions] Let $u\colon M\to\R$ be a temporal function. As $\smash{\vert\rmd u\vert_g^{p-2}\,\rmd u}$ and $g^{-1}$ are locally bounded, \eqref{Eq:Weak} extends to every Lipschitz function $\varphi\colon M\to\R$ with compact support and, for $U\Subset M$ open, to every $\smash{\varphi\in H_0^1(U)}$.
\end{remark}

\begin{definition}[$p$-Harmonic temporal functions]\label{def:pharm} A temporal function $u\colon M\to\R$ is called \textbf{$p$-harmonic} if $\smash{\Box_p u = 0}$ in $\scrD'(M)$.
\end{definition}

Given $x\in M$ and $\zeta\in T_x^*M$ future-directed timelike, we define
\begin{align}\label{Eq:DefA}
    a^\alpha(x,\zeta) := \big\vert\!\det g(x)\big\vert^{1/2}\,\big\vert\zeta\big\vert_{g(x)}^{p-2}\,g^{\alpha\beta}(x)\,\zeta_\beta,
\end{align}
where $\vert\!\det g\vert$ means the modulus of the determinant of $g$. Then \eqref{Eq:Weak} becomes
\begin{align*}
    \big\langle\Box_pu,\varphi\big\rangle = -\int_Ma^\alpha(\cdot,\rmd u)\,\partial_\alpha\varphi\d\Leb^{1+d};
\end{align*}
$p$-harmonicity of a temporal function $u\colon M\to\R$ means $u$ is a weak solution of
\begin{align}\label{Eq:OrigWeakSol}
    \partial_\alpha\big[a^\alpha(\cdot,\rmd u)\big]=0\quad\textnormal{in }M.
\end{align}
Clearly, $\smash{a^\alpha(x,\zeta)}$ depends smoothly on $\zeta\in T_x^*M$ as long as the latter ranges over the future timelike cone in $T_x^*M$; in $x\in M$, it inherits the regularity of $g$, and so do all its derivatives with respect to $\zeta$.

The following ellipticity of $\smash{\Box_p}$ on temporal functions was observed by Beran et al.~\cite{beran-braun-calisti-gigli-mccann-ohanyan-rott-samann+-}; the computation can be found in the proofs of Braun et al.~\cite[Prop.~9]{braun-gigli-mccann-ohanyan-samann+} or Ohanyan \cite{ohanyan2026+}*{Thm.~3.1}, which carry over verbatim to merely continuous $g$.

\begin{lemma}[Ellipticity]\label{lem:ell} Given $x\in M$ and $\zeta\in T_x^*M$ future-directed timelike,
\begin{align*}
    \partial_{\zeta_\beta}a^\alpha(x,\zeta) = \big\vert\!\det g(x)\big\vert^{1/2}\,\big\vert\zeta\big\vert_{g(x)}^{p-4}\,\big[(2-p)\,g^{\alpha\gamma}(x)\,\zeta_\gamma\,g^{\beta\delta}(x)\,\zeta_\delta+\big\vert\zeta\big\vert_{g(x)}^2\,g^{\alpha\beta}(x)\big].
\end{align*}

This matrix is symmetric and positive definite; more precisely, if $\smash{\scrK\subset M\times\R^{1+d}}$ is a compact set such that $\smash{\zeta}$ is future-directed timelike for every $(x,\zeta)\in\scrK$, there exist constants $\lambda,\Lambda > 0$ such that for every $(x,\zeta)\in\scrK$ and every $\smash{\xi\in\R^{1+d}}$,
\begin{align}\label{Eq:EllBounds}
    \lambda\,\vert\xi\vert^2\leq\partial_{\zeta_\beta}a^\alpha(x,\zeta)\,\xi_\alpha\,\xi_\beta\leq\Lambda\,\vert\xi\vert^2.
\end{align}
\end{lemma}

\begin{remark}[Ellipticity and the range of $p$]\label{Re:Whyelliptic} Let $p\in\R$. The identity for $\smash{\partial_\zeta a}$ from \cref{lem:ell} explains our range of $p$. Given $x\in M$, in a $g(x)$-orthonormal frame adapted to $\zeta\in T_x^*M$ future-directed timelike, the matrix $\smash{\partial_\zeta a(x,\zeta)}$ is diagonal with entries $1-p,1,\dots,1$ times the positive factor $\smash{\vert\!\det g(x)\vert^{1/2}\,\vert\zeta\vert_{g(x)}^{p-2}}$. Hence, the linearization of $\smash{\Box_p}$ at a temporal function $u\colon M\to\R$ is elliptic if $p<1$, degenerates in the direction of $\rmd u$ if $p=1$, and is hyperbolic if $p>1$; for $p=2$, it is the wave operator. Compare this to the Riemannian $p$-Laplacian, where the eigenvalues are multiples of $p-1,1,\dots,1$, leading to the elliptic range $p\in(1,\infty)$.
\end{remark}

Uniform ellipticity in the sense of \eqref{Eq:EllBounds} requires covectors to stay in a compact subset of the open cone of future-directed timelike covectors. The following lemma provides such sets for a given temporal function.

\begin{lemma}[Compact convex sets of timelike covectors]\label{lem:K} Let $u\colon M\to\R$ be a temporal function and let $x_0\in M$. Then there exist an open ball $B\Subset M$ around $x_0$ of radius at most $1$ and a compact convex set $\smash{\scrC\subset\R^{1+d}}$ such that 
\begin{enumerate}[label=\textnormal{(\roman*)}]
    \item $\rmd u(x)\in \scrC$ for $\smash{\vol}$-a.e.~$x\in B$ and \item for every $\smash{y\in \overline{B}}$, every element of $\scrC$ constitutes a future-directed timelike covector at $y$.
\end{enumerate}
\end{lemma}

\begin{proof} Let $B_0\Subset M$ be an open ball around $x_0$. Let $L>0$ designate a Lipschitz constant of $u$ on $B_0$ and let $c$ be the constant from \cref{Def:Temporal} relative to $\smash{\overline{B}_0}$. By uniform continuity of $\smash{g^{-1}}$ and the time orientation $X$ on $\smash{\overline{B}_0}$, there is an open ball $B\subset B_0$ around $x_0$ of radius at most $1$ such that for every $\smash{x,z\in\overline{B}}$,
\begin{align}\label{Eq:Bds}
\begin{split}
    \big\vert g^{-1}(z)-g^{-1}(x)\big\vert &\leq \frac{c}{2L^2},\\
    \big\vert X(z)-X(x)\big\vert &\leq \frac{c}{2L}.
\end{split}
\end{align}

We now verify the claim for
\begin{align*}
    \scrC := \big\{\zeta\in \overline{B}_L(0): g^{-1}(y)(\zeta,\zeta)\leq -\delta \textnormal{ and }\zeta(X(y))\geq \delta \textnormal{ for every }y\in\overline{B}\big\},
\end{align*}
where $\smash{\delta:=c/2}$. Clearly, $\scrC$ is compact, as it is closed and bounded in $\smash{\R^{1+d}}$. By definition, every element of $\scrC$ is a future-directed timelike covector at every $\smash{y\in\overline{B}}$. Moreover, $\scrC$ is convex: given $\smash{y\in\overline{B}}$, the set of all future-directed timelike covectors $\smash{\zeta\in T_y^*M}$ with $\smash{\vert \zeta\vert_{g(y)}\geq\delta^{1/2}}$ is convex by the reverse Minkowski inequality and $\scrC$ is the intersection of these sets with half-spaces and a ball. Finally, let $x\in B$ be a point at which $u$ is differentiable (which is the case $\vol$-a.e.~by Rademacher's theorem, since $u$ is locally Lipschitz) and where \eqref{Eq:LowBo} holds. Then $\smash{\rmd u(x)\in \overline{B}_L(0)}$. The remaining estimates are direct consequences of this inequality together with \eqref{Eq:Bds} and \eqref{Eq:LowBo}. This establishes $\rmd u(x)\in\scrC$.
\end{proof}

\section{Some auxiliary estimates}\label{Sec:Lemmas}

Throughout the rest of this paper, $g$ will be a Lorentzian metric on $M$ of class $\smash{C^{0,\alpha}_\loc}$, where $\alpha\in(0,1)$. We retain the notation from \cref{Sec:Prelim}.

\subsection{Truncation and linearization} Let us consider the coefficient $a$ from \eqref{Eq:DefA}. In this section, we will truncate it to remove its degeneracy at the light cone: we restrict it to a compact convex set of future-directed timelike covectors and extend it from there to a globally Lipschitz and strongly monotone map on $\smash{\R^{1+d}}$. This will allow us to invoke general existence, uniqueness, and regularity theory for elliptic PDEs with strongly monotone coefficients, both in \cref{Sec:Reg} and the construction of $p$-harmonic temporal functions in \cref{Sec:Exist}. We adopt the framework of \cref{lem:K}: we fix an open ball $B\Subset M$ and let $\smash{\scrC\subset\R^{1+d}}$ be compact and convex such that for every $\smash{x\in\overline{B}}$, every element of $\scrC$ is a future-directed timelike covector at $x$. 

\begin{lemma}[Local Lipschitz continuity and strong monotonicity]\label{Le:LocLip} There exist constants $c>0$ and $\Lambda_0\geq 4c$ such that for every $\smash{x\in\overline{B}}$ and every $\zeta,\zeta'\in\scrC$,
\begin{align}\label{Eq:OnK}
\begin{split}
    \big[a^\alpha(x,\zeta)-a^\alpha(x,\zeta')\big]\,(\zeta-\zeta')_\alpha&\geq 4c\,\vert \zeta-\zeta'\vert^2,\\
    \vert a(x,\zeta)-a(x,\zeta')\vert&\leq\Lambda_0\,\vert \zeta-\zeta'\vert.
\end{split}
\end{align}
\end{lemma}

\begin{proof} The set $\smash{\scrK:=\overline{B}\times\scrC}$ is compact in $\smash{\R^{1+d}\times \R^{1+d}}$. Given $\smash{x\in\overline{B}}$ and $\zeta,\zeta'\in\scrC$, the segment from $\zeta'$ to $\zeta$ lies in $\scrC$ by convexity, hence
\begin{align*}
    a^\alpha(x,\zeta)-a^\alpha(x,\zeta') = \int_{[0,1]}\partial_{\zeta_\beta}a^\alpha(x,\zeta'+t\,(\zeta-\zeta'))\,(\zeta-\zeta')_\beta\d t.
\end{align*}
The lower bound of \eqref{Eq:EllBounds} gives the first inequality from \eqref{Eq:OnK} with an appropriate choice of $c$. Letting $\Lambda_0>0$ be the maximum of the norm of $\partial_\zeta a$ on $\scrK$, which clearly satisfies $\Lambda_0\geq 4c$, the second estimate from \eqref{Eq:OnK} follows readily.
\end{proof}

Fix constants $c$ and $\Lambda_0$ as above and abbreviate $\lambda := \Lambda_0^2/4c$. Furthermore, let $\smash{P\colon \R^{1+d}\to\scrC}$ denote the nearest-point projection onto $\scrC$. Define the \textbf{truncation} $\smash{\overline{a}\colon \overline{B}\times \R^{1+d}\to\R^{1+d}}$ of $a$ by
\begin{align*}
    \overline{a}(x,\zeta) := a(x,P\zeta) + \lambda\,(\zeta-P\zeta).
\end{align*}

\begin{lemma}[Global Lipschitz and Hölder continuity and strong monotonicity]\label{lem:trunc} The following properties hold for every $\smash{x,x'\in\overline{B}}$ and every $\smash{\zeta,\zeta'\in\R^{1+d}}$.
\begin{enumerate}[label=\textnormal{\textcolor{black}{(}\roman*\textcolor{black}{)}}]
    \item\label{La:T1} If $\zeta\in\scrC$, we have $\overline{a}(x,\zeta) = a(x,\zeta)$.
    \item\label{La:T2} Setting $\Lambda := \Lambda_0+\lambda$, we have the Lipschitz estimate
    \begin{align*}
        \big\vert\overline{a}(x,\zeta) -\overline{a}(x,\zeta')\big\vert \leq \Lambda\,\vert\zeta-\zeta'\vert.
    \end{align*}
    \item\label{La:T3} We have the strong monotonicity estimate
    \begin{align*}
        \big[\overline{a}^\alpha(x,\zeta)-\overline{a}^\alpha(x,\zeta')\big]\,(\zeta-\zeta')_\alpha\geq c\,\vert \zeta-\zeta'\vert^2.    
    \end{align*}
     \item\label{La:T4} There exists a constant $K>0$ depending only on $\scrC$, the norms $\smash{\Vert g\Vert_{C^{0,\alpha}(\overline{B})}}$ and $\smash{\Vert g^{-1}\Vert_{C^0(\overline{B})}}$, and the positive number
\begin{align*}
    \delta := \min\{{-g^{-1}(x)(w,w)} : x\in\overline{B},\, w\in\scrC\}
\end{align*}
such that 
\begin{align*}
    \big\vert\overline{a}(x,\zeta)-\overline{a}(x',\zeta)\big\vert\leq K\,\vert x-x'\vert^\alpha.
\end{align*}
\end{enumerate}
\end{lemma}

\begin{proof} \ref{La:T1} This is clear. 

\ref{La:T2} follows from \cref{Le:LocLip} and $1$-Lipschitz continuity of $P$ and $\Id-P$.

\ref{La:T4} Given $\smash{\zeta\in\R^{1+d}}$, note that $\smash{\overline{a}(x,\zeta)}$ depends on $\smash{x\in \overline{B}}$ only through the term $a(x,w)$, where $w:=P\zeta$. By \eqref{Eq:DefA}, we have $a(x,w) = A(g(x),w)$, where
\begin{align*}
    A(G,w) := \vert\!\det G\vert^{1/2}\,\big\vert w\big\vert_{G}^{p-2}\,G^{-1}w
\end{align*}
is smooth on the set of all pairs $(G,w)$ such that $\smash{G\in\R^{(1+d)\times(1+d)}}$ is a symmetric matrix of Lorentzian signature and $\smash{G^{-1}(w,w)<0}$, which is open in the space of symmetric matrices times $\smash{\R^{1+d}}$. Let $L>0$ denote the maximum of $\smash{\Vert g\Vert_{C^0(\overline{B})}}$ and $\smash{\Vert g^{-1}\Vert_{C^0(\overline{B})}}$. Let $S$ denote the set of all symmetric matrices $\smash{G\in\R^{(1+d)\times(1+d)}}$ of Lorentzian signature with $\vert G\vert\leq L$, $\smash{\vert G^{-1}\vert\leq L}$, and $\smash{G^{-1}(w,w)\leq-\delta}$ for every $w\in\scrC$. Then $S\times\scrC$ is a compact subset of the domain of $A$ that contains $\smash{g(\overline{B})\times\scrC}$. Thus, $A$ is Lipschitz continuous on $S\times\scrC$ with a constant depending only on $L$, $\delta$, and $\scrC$. Together with the $\alpha$-Hölder continuity of $g$ on $\smash{\overline{B}}$, this gives the desired claim.

\ref{La:T3} We abbreviate
\begin{align*}
\Delta P &:= P\zeta-P\zeta',\\
\Delta (\Id-P) &:= (\zeta-P\zeta)-(\zeta'-P\zeta'),    
\end{align*}
so that $\zeta-\zeta' = \Delta P+\Delta (\Id-P)$. Recall the operator $\mathrm{Id}-P$ is firmly nonexpansive, cf.~e.g.~Bauschke--Combettes \cite[Prop.~4.16]{bauschke-combettes2017}, which means
\begin{align*}
    \delta^{\alpha\beta}\,\Delta(\Id-P)_\beta\,(\zeta-\zeta')_\alpha\geq \vert\Delta(\Id-P)\vert^2.
\end{align*}
Combining this with \cref{Le:LocLip} and Young's inequality,
\begin{align*}
    &\big[\overline{a}^\alpha(x,\zeta)-\overline{a}^\alpha(x,\zeta')\big]\,(\zeta-\zeta')_\alpha\\
    &\qquad\qquad= \big[a^\alpha(x,P\zeta)-a^\alpha(x,P\zeta')\big]\,\Delta P_\alpha\\
    &\qquad\qquad\qquad\qquad + \big[a^\alpha(x,P\zeta)-a^\alpha(x,P\zeta')\big]\,\Delta(\Id- P)_\alpha\\
    &\qquad\qquad\qquad\qquad+\lambda\,\delta^{\alpha\beta}\,\Delta(\Id-P)_\beta\,(\zeta-\zeta')_\alpha\\
    &\qquad\qquad\geq 4c\,\vert\Delta P\vert^2-\Lambda_0\,\vert\Delta P\vert\,\vert\Delta (\Id-P)\vert+\lambda\,\vert\Delta (\Id-P)\vert^2\\
    &\qquad\qquad\geq 2c\,\vert\Delta P\vert^2+\Big[\lambda-\frac{\Lambda_0^2}{8c}\Big]\,\vert\Delta (\Id-P)\vert^2\\
    &\qquad\qquad\geq 2c\,\big[\vert\Delta P\vert^2+\vert\Delta (\Id-P)\vert^2\big]\\
    &\qquad\qquad\geq c\,\vert \zeta-\zeta'\vert^2;
\end{align*}
in the second last step, we used $\lambda-\Lambda_0^2/8c = \Lambda_0^2/8c\geq 2c$ by our choice of $\lambda$.
\end{proof}

The rest of this section will be dedicated to establishing auxiliary estimates for coefficients satisfying the estimates from \cref{lem:trunc} and for weak solutions of the associated elliptic PDEs.

The truncated map $\overline{a}$ is merely Lipschitz in its second argument. The following simple observation allows us to linearize differences of such maps nevertheless.

\begin{lemma}[Linearization]\label{lem:lin} Let $c>0$ and $\Lambda \geq c$. Assume $\smash{F\colon \R^{1+d}\to\R^{1+d}}$ is a $\Lambda$-Lipschitz map such that for every $\smash{\zeta,\zeta'\in\R^{1+d}}$,
\begin{align*}
    \big[F^\alpha(\zeta)-F^\alpha(\zeta')\big]\,(\zeta-\zeta')_\alpha\geq c\,\vert\zeta-\zeta'\vert^2.
\end{align*}
Define $\sfa\colon \R^{1+d}\times\R^{1+d}\to \R^{(1+d)\times (1+d)}$ by
\begin{align*}
\sfa(\zeta,\zeta')^{\alpha\beta}:= \frac{\Lambda^2}{c}\,\big[\delta^{\alpha\beta} - e(\zeta,\zeta')^\alpha\,e(\zeta,\zeta')^\beta\big] + \frac{[F^\alpha(\zeta') - F^\alpha(\zeta)]\,e(\zeta,\zeta')^\beta}{\vert\zeta'-\zeta\vert}
\end{align*}
with the convention $0/0 := 0$, where $e\colon \R^{1+d}\times \R^{1+d}\to\S^d\cup\{0\}$ is given by
\begin{align*}
    e(\zeta,\zeta'):= \begin{cases}
        \dfrac{\zeta'-\zeta}{\vert\zeta'-\zeta\vert} & \textnormal{\textit{if }}\zeta\neq\zeta',\\
        0 & \textnormal{\textit{otherwise}}.
    \end{cases}
\end{align*}
Then for every $\smash{\zeta,\zeta',\xi\in\R^{1+d}}$,
\begin{align*}
    \sfa(\zeta,\zeta')^{\alpha\beta}\,(\zeta'-\zeta)_\beta &= F^\alpha(\zeta')-F^\alpha(\zeta),\\
    \big\vert\sfa(\zeta,\zeta')\big\vert&\leq(1+d)^{1/2}\, \frac{\Lambda^2}{c}+\Lambda,\\
    \sfa(\zeta,\zeta')^{\alpha\beta}\,\xi_\alpha\,\xi_\beta&\geq \frac{c}{2}\,\vert\xi\vert^2.
\end{align*}
\end{lemma}

\begin{proof} All claims are clear on the diagonal, so assume $\zeta\neq\zeta'$ throughout the proof.

The first summand of $\sfa(\zeta,\zeta')$ is a multiple of the orthogonal projection onto the orthogonal complement of $e(\zeta,\zeta')$, thus vanishes at $\smash{\zeta'-\zeta}$. The second summand maps $\zeta'-\zeta$ to $F(\zeta')-F(\zeta)$. This establishes the first claim.

The second claim follows since the orthogonal projection onto a $d$-dimensional subspace has norm $\smash{d^{1/2}}$ and since the norm of the second summand is at most $\Lambda$ by the hypothesized $\Lambda$-Lipschitz continuity of $F$.

For the third claim, we write $\xi = s\,e(\zeta,\zeta') + \xi_\perp$, where $s\in\R$ and $\smash{\xi_\perp\in \R^{1+d}}$ is orthogonal to $e(\zeta,\zeta')$. Then
\begin{align*}
    e(\zeta,\zeta')^\beta\,\xi_\beta = s\,e(\zeta,\zeta')^\beta\,e(\zeta,\zeta')_\beta + e(\zeta,\zeta')^\beta\,\xi_{\perp,\beta} =s
\end{align*}
and consequently
\begin{align*}
    \sfa(\zeta,\zeta')^{\alpha\beta}\,\xi_\alpha\,\xi_\beta &= \frac{\Lambda^2}{c}\,\vert\xi_\perp\vert^2 + \frac{s^2\,[F^\alpha(\zeta') - F^\alpha(\zeta)]\,e(\zeta,\zeta')_\alpha}{\vert\zeta-\zeta'\vert}\\
    &\qquad\qquad + \frac{s\,[F^\alpha(\zeta')-F^\alpha(\zeta)]\,\xi_{\perp,\alpha}}{\vert \zeta-\zeta'\vert}\\
    &\geq \frac{\Lambda^2}{c}\,\vert\xi_\perp\vert^2 + c\,s^2 -\Lambda\,\vert s\vert\,\vert\xi_\perp\vert\\
    &\geq \frac{c}{2}\,s^2 + \frac{\Lambda^2}{2c}\,\vert\xi_\perp\vert^2\\
    &\geq \frac{c}{2}\,\vert\xi\vert^2;
\end{align*}
here, we used Young's inequality in the second last step and $\Lambda \geq c$ in the last.
\end{proof}

\subsection{Iteration and Campanato's criterion} The following result is a special case of Giaquinta--Martinazzi \cite{giaquinta-martinazzi2012}*{Lem.~5.13} with $\varepsilon$ replaced by $0$.

\begin{lemma}[Iteration]\label{lem:iter} Let $R_0>0$ and let $\Theta\colon(0,R_0]\to\R_+$ be nondecreasing. Furthermore, let $A_0\geq 1$, $B_0\in \R_+$, and $a,b>0$ with $a<b$ such that for every $\varrho,R\in (0,R_0]$ such that $\varrho\leq R$,
\begin{align}\label{Eq:Hyyyyyy}
    \Theta(\varrho)\leq A_0\,\Big[\frac{\varrho}{R}\Big]^b\,\Theta(R)+B_0\,R^a.
\end{align}
Then there exists a constant $C>0$ depending only on $A_0$, $a$, and $b$ such that for every $\varrho,R\in(0,R_0]$ with $\varrho\leq R$,
\begin{align*}
    \Theta(\varrho)\leq C\,\Big[\Big[\frac{\varrho}{R}\Big]^a\,\Theta(R)+B_0\,\varrho^a\Big].
\end{align*}
\end{lemma}

Given $x\in \R^{1+d}$, $\varrho > 0$, and $f\in L^2(B_\varrho(x))$, the \textbf{mean} of $f$ on $B_\varrho(x)$ is
\begin{align}\label{Eq:MEAN!}
\begin{split}
    f_{x,\varrho} &:= \frac{1}{\Leb^{1+d}[B_\varrho(x)]}\int_{B_\varrho(x)} f\d\Leb^{1+d}\\
    &\phantom{:}=\argmin\!\Big\{\!\int_{B_\varrho(x)} \vert f-k\vert^2\d\Leb^{1+d}: k\in\R^{1+d}\Big\}.
\end{split}
\end{align}

We will use the following two classical facts. The first is Campanato's integral characterization of Hölder continuity \cite{campanato1963}, see also Meyers \cite{meyers1964}. The second is a decay estimate for solutions of linear elliptic equations with constant coefficients. We refer to Giaquinta--Martinazzi \cite[§5]{giaquinta-martinazzi2012} for relevant background.

\begin{lemma}[Campanato's criterion for Hölder continuity]\label{lem:campanato} Let $\smash{z\in\R^{1+d}}$, $s>0$, $\sigma\in(0,1)$, and $H\in\R_+$. Assume $\smash{f\in L^2(B_{2s}(z);\R^{1+d})}$ satisfies
\begin{align*}
    \int_{B_\varrho(y)}\big\vert f-f_{y,\varrho}\big\vert^2\d\Leb^{1+d}\leq H^2\,\varrho^{1+d+2\sigma}
\end{align*}
for every $y\in B_s(z)$ and every $\varrho\in(0,s]$. Then $f$ has an $\smash{\Leb^{1+d}}$-representative which belongs to $\smash{C^{0,\sigma}(\overline{B}_s(z))}$ with $\sigma$-Hölder seminorm at most $CH$, where $C>0$ is a constant depending only on $d$ and $\sigma$.
\end{lemma}

\begin{proof} Both the hypothesis and the conclusion are invariant under replacing $f$ by $f(z+s\,\cdot)$, which is in $\smash{L^2(B_2(0),\R^{1+d})}$, and $H$ by $\smash{s^\sigma H}$. Thus, by translation and dilation, it suffices to consider $z=0$ and $s=1$.

Let $y\in B_1(0)$ and $\varrho>0$. Let $\smash{k_{y,\varrho}\in\R^{1+d}}$ be the mean of $f$ on $\smash{B_1(0)\cap B_\varrho(y)}$. Recall $k_{y,\varrho}$ minimizes the distance to constants in $\smash{L^2(B_1(0)\cap B_\varrho(y))}$. If $\varrho \leq 1$, we choose $\smash{f_{y,\varrho}}$ as a competitor in \eqref{Eq:MEAN!} and obtain
\begin{align*}
    \int_{B_1(0)\cap B_\varrho(y)}\big\vert f-k_{y,\varrho}\big\vert^2\d\Leb^{1+d}\leq\int_{B_\varrho(y)}\big\vert f-f_{y,\varrho}\big\vert^2\d\Leb^{1+d}\leq H^2\,\varrho^{1+d+2\sigma}.
\end{align*}
If $\varrho>1$, we choose $\smash{f_{0,1}}$ as a competitor in \eqref{Eq:MEAN!} to get
\begin{align*}
    \int_{B_1(0)\cap B_\varrho(y)}\big\vert f-k_{y,\varrho}\big\vert^2\d\Leb^{1+d}\leq\int_{B_1(0)}\big\vert f-f_{0,1}\big\vert^2\d\Leb^{1+d}\leq H^2\leq H^2\,\varrho^{1+d+2\sigma}.
\end{align*}
This shows the Campanato seminorm
\begin{align*}
    \sup\!\Big\{\Big[\varrho^{-1-d-2\sigma}\int_{B_1(0)\cap B_\varrho
    (y)} \big\vert f-k_{y,\varrho}\big\vert^2\d\Leb^{1+d}\Big]^{1/2}: y \in B_1(0)\textnormal{ and }\varrho>0\Big\}
\end{align*}
of $f$ on $B_1(0)$ with exponents $2$ and $1+d+2\sigma$ is at most $H$. By Campanato's theorem, cf.~Giaquinta--Martinazzi \cite[Thm.~5.5]{giaquinta-martinazzi2012} (applied to $n$, $\lambda$, and $p$ replaced by $1+d$, $1+d+2\sigma$, and $2$, respectively), $f$ has a representative that is $\sigma$-Hölder continuous on $\smash{\overline{B}_1(0)}$ with seminorm bounded from above by $CH$, where the constant $C>0$ depends only on $d$ and $\sigma$.
\end{proof}

\begin{lemma}[Decay for elliptic PDEs with constant coefficients]\label{lem:const} Fix a matrix $\smash{\sfb \in\R^{(1+d)\times(1+d)}}$ and constants $\lambda,\Lambda >0$ such that $\vert\sfb\vert\leq \Lambda$ and $\smash{\sfb^{\alpha\beta}\,\xi_\alpha\xi_\beta\geq \lambda\,\vert\xi\vert^2}$ for every $\smash{\xi\in\R^{1+d}}$. Then there is a constant $C\geq 1$ depending only on $d$, $\lambda$, and $\Lambda$ such that for every $\smash{y\in\R^{1+d}}$ and every $R>0$, every weak solution $\smash{v\in H^1(B_R(y))}$ of the elliptic equation
\begin{align*}
    \partial_\alpha\big[\sfb^{\alpha\beta}\,\partial_\beta v\big]=0\quad\textnormal{\textit{in} }B_R(y)
\end{align*}
satisfies, for every $\varrho\in(0,R]$,
\begin{align*}
 \big\Vert\rmd v-(\rmd v)_{y,\varrho}\big\Vert_{L^2(B_\varrho(y))}^2\leq C\,\Big[\frac{\varrho}{R}\Big]^{3+d}\big\Vert\rmd v-(\rmd v)_{y,R}\big\Vert_{L^2(B_R(y))}^2.
\end{align*}
\end{lemma}

\begin{proof} Since the coefficients are constant, $v$ belongs to $\smash{H^2_\loc(B_R(y))}$ by the method of difference quotients, cf.~e.g.~Giaquinta--Martinazzi \cite[Thm.~4.9]{giaquinta-martinazzi2012}. Testing the equation with partial derivatives of test functions and integrating by parts, we see the partial derivatives of $v$ are weak solutions of the same equation.

Let $\varrho\in(0,R)$ and $R'\in(\varrho,R)$. Then \cite[Prop.~5.8]{giaquinta-martinazzi2012} applied to every partial derivative of $v$  yields
\begin{align*}
\big\Vert\rmd v- (\rmd v)_{y,\varrho}\big\Vert_{L^2(B_\varrho(y))}^2  &\leq C\,\Big[\frac{\varrho}{R'}\Big]^{3+d}\big\Vert\rmd v-(\rmd v)_{y,R'}\big\Vert^2_{L^2(B_{R'}(y))},
\end{align*}
where $C\geq 1$ is a constant depending only on $d$, $\lambda$, and $\Lambda$. Sending $R'\to R$ yields the claim for $\varrho < R$; the case $\varrho = R$ is trivial.
\end{proof}

\subsection{A Hölder gradient estimate} Now we derive a gradient estimate for weak solutions $\smash{u\in H^1(B_1(0))}$ of the PDE
\begin{align}\label{Eq:TruncUnit}
    \partial_\alpha\big[b^\alpha(\cdot,\rmd u)\big]=0\quad\textnormal{in }B_1(0),
\end{align}
where $\smash{b\colon B_1(0)\times \R^{1+d}\to\R^{1+d}}$ satisfies estimates analogous to those of $\smash{\overline{a}}$ from \cref{lem:trunc}. The corresponding next lemma is our main technical ingredient. The proof follows the scheme of Giaquinta--Martinazzi \cite[Thm.~8.32]{giaquinta-martinazzi2012}; the differences are that our equation is not variational and that $b$ is merely Lipschitz in its second argument, which is handled by \cref{lem:lin}.

\begin{lemma}[Hölder gradient estimate]\label{lem:camp} Let $c>0$, $\Lambda\geq c$, $K\in \R_+$, and $\alpha\in(0,1)$. Let $\smash{b\colon B_1(0)\times\R^{1+d}\to\R^{1+d}}$ be a map satisfying the following estimates for every $x,x'\in B_1(0)$ and every $\smash{\zeta,\zeta'\in\R^{1+d}}$:
\begin{align}\label{Eq:Struct}
\begin{split}
    \big\vert b(x,\zeta)-b(x,\zeta')\big\vert&\leq\Lambda\,\vert \zeta-\zeta'\vert,\\
    \big[b^\alpha(x,\zeta)-b^\alpha(x,\zeta')\big]\,(\zeta-\zeta')_\alpha&\geq c\,\vert \zeta-\zeta'\vert^2,\\
    \big\vert b(x,\zeta)-b(x',\zeta)\big\vert&\leq K\,\vert x-x'\vert^\alpha.
\end{split}
\end{align}
Then there exist $\beta\in(0,\alpha]$ and a constant $C>0$ depending only on $d$, $c$, $\Lambda$, and $\alpha$ with the following property. For every weak solution $\smash{u\in H^1(B_1(0))}$ of \eqref{Eq:TruncUnit}, the differential $\rmd u$ has a $\beta$-Hölder continuous representative on $\smash{\overline{B}_{1/2}(0)}$ which satisfies, for every $\smash{\zeta_0\in\R^{1+d}}$,
\begin{align*}
   \big\Vert \rmd u-\zeta_0\big\Vert_{C^{0,\beta}(\overline{B}_{1/2}(0))} \leq C\,\big[\big\Vert\rmd u-\zeta_0\big\Vert_{L^2(B_1(0))}+K\big].
\end{align*}
\end{lemma}

\begin{proof} The proof proceeds in five steps. First, we recall a decay estimate for weak solutions of linear elliptic equations with merely measurable coefficients, which follows from the theory of De Giorgi--Nash--Moser. Second, we freeze the first argument of $b$ at a point; by \cref{lem:lin}, difference quotients of weak solutions of the frozen equation solve linear equations of the former type, which transfers the decay estimate to their differentials. Third, on every small ball we compare $u$ to the solution of the frozen equation with the same boundary values; the error is controlled by the Hölder continuity of $b$ in its first argument. Finally, we combine the last two steps by means of \cref{lem:iter} and conclude using \cref{lem:campanato}.

Throughout, $C\geq 1$ denotes a constant depending only on $d$, $c$, $\Lambda$, and $\alpha$ whose value may change from line to line.

\emph{Step 1. Linear equations with measurable coefficients.} Let $\smash{y\in \R^{1+d}}$ and $R>0$. Let $\smash{\sfb\colon B_R(y)\to\R^{(1+d)\times(1+d)}}$ be a Lebesgue measurable map satisfying the last two inequalities from \cref{lem:lin}; that is, for every $\smash{\xi\in\R^{1+d}}$, 
\begin{align*}
\vert\sfb\vert &\leq(1+d)^{1/2}\,\frac{\Lambda^2}{c}+\Lambda,\\
\sfb^{\alpha\beta}\,\xi_\alpha\xi_\beta &\geq\frac{c}{2}\,\vert\xi\vert^2.    
\end{align*}
We claim there exists $\gamma\in(0,1)$ depending only on $d$, $c$, and $\Lambda$ such that every weak solution $\smash{z\in H^1(B_R(y))}$ of
\begin{align}\label{Eq::::PDE}
    \partial_\alpha\big[\sfb^{\alpha\beta}\,\partial_\beta z\big]=0\quad\textnormal{in }B_R(y)
\end{align}
satisfies, for every $\varrho\in(0,R]$ and every $k\in\R$,
\begin{align}\label{Eq:DGNM}
\big\Vert z-z_{y,\varrho}\big\Vert^2_{L^2(B_\varrho(y))}\leq C\,\Big[\frac{\varrho}{R}\Big]^{1+d+2\gamma}\big\Vert z-k\big\Vert^2_{L^2(B_R(y))}.
\end{align}
Indeed, the function $z-k$ of regularity $H^1$ is a weak solution of the same equation \eqref{Eq::::PDE}. We apply two results from the theory of De Giorgi--Nash--Moser; in both, the constants depend only on $d$, $c$, and $\Lambda$, and not on $R$, since the equation has no lower order terms and is thus invariant under rescaling. First, the local maximum principle, cf.~e.g.~Gilbarg--Trudinger \cite[Thm.~8.17]{gilbarg-trudinger2001}, applied on $B_R(y)$ with exponent $2$ and radius $R/2$ to the subsolutions $z-k$ and $-z+k$, yields
\begin{align*}
    \sup\{\pm(z(x)-k) : x\in B_{R/2}(y)\}\leq C\,R^{-(1+d)/2}\,\big\Vert (z-k)^\pm\big\Vert_{L^2(B_R(y))}
\end{align*}
and consequently
\begin{align*}
   \big\Vert z-k\big\Vert_{C^0(B_{R/2}(y))}^2\leq C\,R^{-1-d}\int_{B_R(y)}\vert z-k\vert^2\d\Leb^{1+d}.
\end{align*}
Second, the interior Hölder estimate \cite[Thm.~8.22]{gilbarg-trudinger2001}, applied to $z-k$ on the ball $B_{R/2}(y)$, yields a constant $\gamma\in(0,1)$ such that for every $\varrho\in(0,R/2]$,
\begin{align*}
    \mathrm{osc}_{B_\varrho(y)}\,z = \mathrm{osc}_{B_\varrho(y)}\,(z-k)\leq C\,\Big[\frac{\varrho}{R}\Big]^\gamma\Vert z-k\Vert_{C^0(B_{R/2}(y))}.
\end{align*}
Given $\varrho\in(0,R/2]$, combining these two estimates implies \eqref{Eq:DGNM} since the left-hand side of \eqref{Eq:DGNM} is at most $\smash{\Leb^{1+d}[B_\varrho(y)]}$ times the square of the oscillation of $z$ on $B_\varrho(y)$. For $\varrho\in(R/2,R]$, the estimate \eqref{Eq:DGNM} holds as soon as $\smash{C\geq 2^{3+d}}$, by the minimality property \eqref{Eq:MEAN!} of the mean.

\emph{Step 2. Decay for the frozen equation.} We fix $y_1\in B_{1/2}(0)$ and $R\in(0,1/2]$. Let $\smash{v\in H^1(B_R(y_1))}$ be a weak solution of
\begin{align}\label{Eq:Frozen}
    \partial_\alpha\big[b^\alpha(y_1,\rmd v)\big]=0\quad\textnormal{in }B_R(y_1).
\end{align}
We claim that for every $\varrho\in(0,R]$,
\begin{align}\label{Eq:Decay}
\begin{split}
    \big\Vert\rmd v-(\rmd v)_{y_1,\varrho}\big\Vert^2_{L^2(B_\varrho(y_1))} \leq C\,\Big[\frac{\varrho}{R}\Big]^{1+d+2\gamma}\big\Vert\rmd v-(\rmd v)_{y_1,R}\big\Vert^2_{L^2(B_R(y_1))}.
    \end{split}
\end{align}
Let $R'\in(0,R)$, let $\smash{e_\kappa\in\R^{1+d}}$ be a coordinate direction, and let $h\in (R'-R,R-R')\setminus\{0\}$. Consider the difference quotient $\smash{z_h\in H^1(B_{R'}(y_1))}$ with
\begin{align*}
    z_h := \frac{v(\cdot+he_\kappa)-v}{h}.
\end{align*}
Given $\smash{\varphi\in C_\comp^\infty(B_{R'}(y_1))}$, testing \eqref{Eq:Frozen} with $\varphi(\cdot-he_\kappa)$ and with $\varphi$ and subtracting, 
\begin{align}\label{Eq:PDE!}
    \int_{B_{R'}(y_1)}\big[b^\alpha(y_1,\rmd v(\cdot +he_\kappa))-b^\alpha(y_1,\rmd v)\big]\,\partial_\alpha\varphi\d\Leb^{1+d} = 0.
\end{align}
Thanks to \eqref{Eq:Struct}, the map $\smash{F\colon\R^{1+d}\to\R^{1+d}}$ defined by $F(\zeta) := b(y_1,\zeta)$ satisfies the hypotheses of \cref{lem:lin}. Let $\sfa$ denote the corresponding matrix-valued map from that lemma and abbreviate $\smash{\sfb_h(y) := \sfa(\rmd v(y),\rmd v(y+he_\kappa))}$, 
a quantity which depends Lebesgue measurably on $\smash{y\in B_{R'}(y_1)}$. As $\rmd v(y+he_\kappa)-\rmd v(y) = h\,\rmd z_h(y)$, \cref{lem:lin} turns the identity \eqref{Eq:PDE!} into
\begin{align*}
    \int_{B_{R'}(y_1)}\sfb_h^{\alpha\beta}\,\partial_\beta z_h\,\partial_\alpha\varphi\d\Leb^{1+d} = 0.
\end{align*}
By the arbitrariness of $\varphi$, this shows $z_h$ is a weak solution of a linear equation of the form \eqref{Eq::::PDE} on $B_{R'}(y_1)$. By \eqref{Eq:DGNM}, for every $\varrho\in(0,R']$ and every $k\in\R$,
\begin{align*}
    \big\Vert z_h-(z_h)_{y_1,\varrho}\big\Vert_{L^2(B_\varrho(y_1))}^2\leq C\,\Big[\frac{\varrho}{R'}\Big]^{1+d+2\gamma}\big\Vert z_h-k\big\Vert_{L^2(B_{R'}(y_1))}^2.
\end{align*}
As $h\to 0$, we know that $z_h\to\partial_\kappa v$ in $\smash{L^2(B_{R'}(y_1))}$, cf.~e.g.~Giaquinta--Martinazzi \cite{giaquinta-martinazzi2012}*{Prop.~4.8}. In particular, the previous estimate passes to the limit. Letting $R'\to R$ subsequently, replacing $k$ by $\smash{(\partial_\kappa v)_{y_1,R}}$, and summing over $\kappa$ yields \eqref{Eq:Decay} for $\varrho\in(0,R)$; for $\varrho=R$, there is nothing to prove.

\emph{Step 3. Comparison.} Let $u$ and $\zeta_0$ be as in the statement. Let $y_1\in B_{1/2}(0)$ and $R\in(0,1/2]$, which guarantees in particular $B_R(y_1)\subset B_1(0)$. By \eqref{Eq:Struct} and the Poincaré inequality on $B_R(y_1)$, the assignment
\begin{align*}
    \langle T\chi,\varphi\rangle := \int_{B_R(y_1)}b^\alpha(y_1,\rmd u+\rmd\chi)\,\partial_\alpha\varphi\d\Leb^{1+d}
\end{align*}
for every $\smash{\varphi\in H_0^1(B_R(y_1))}$ does define a Lipschitz and strongly monotone operator $\smash{T\colon H^1_0(B_R(y_1))\to H^{-1}(B_R(y_1))}$. By Zarantonello's theorem \cite[Thm.~25.B]{zeidler1990b}, a special case of the Browder--Minty theorem \cite{minty1962,browder1963} on existence and uniqueness of solutions to equations with monotone operators, there is a unique $\smash{\chi\in H^1_0(B_R(y_1))}$ such that $T\chi=0$. In other words, $\smash{v\in H^1(B_R(y_1))}$ given by $v := u+\chi$ is a weak solution of \eqref{Eq:Frozen} with Dirichlet boundary conditions $\smash{v-u\in H^1_0(B_R(y_1))}$. 

By \eqref{Eq:Struct}, the map $b$ has at most linear growth in its second argument, uniformly in the first. Hence, by density, the weak formulations of \eqref{Eq:TruncUnit} and \eqref{Eq:Frozen} extend from smooth compactly supported test functions to test functions in $\smash{H^1_0(B_1(0))}$ and $\smash{H^1_0(B_R(y_1))}$, respectively. Thus, testing \eqref{Eq:Frozen} with $u-v$ and \eqref{Eq:TruncUnit} with its extension by zero to $B_1(0)$, which belongs to $\smash{H^1_0(B_1(0))}$, gives
\begin{align*}
    \int_{B_R(y_1)} b^\alpha(y_1,\rmd v)\,\partial_\alpha(u-v)\d\Leb^{1+d} = 0 = \int_{B_R(y_1)} b^\alpha(\cdot,\rmd u)\,\partial_\alpha(u-v)\d\Leb^{1+d}.
\end{align*}
This combines with \eqref{Eq:Struct} to yield
\begin{align*}
    c\,\big\Vert\rmd u-\rmd v\big\Vert_{L^2(B_R(y_1))}^2 &\leq\int_{B_R(y_1)}\big[b^\alpha(y_1,\rmd u)-b^\alpha(y_1,\rmd v)\big]\,\partial_\alpha(u-v)\d\Leb^{1+d}\\
    &= \int_{B_R(y_1)}\big[b^\alpha(y_1,\rmd u)-b^\alpha(\cdot,\rmd u)\big]\,\partial_\alpha(u-v)\d\Leb^{1+d}\\
    &\leq K\,R^\alpha\int_{B_R(y_1)}\vert\rmd u-\rmd v\vert\d\Leb^{1+d}.
\end{align*}
By the Cauchy--Schwarz inequality, this implies
\begin{align}\label{Eq:Comparison}
    \big\Vert\rmd u-\rmd v\big\Vert_{L^2(B_R(y_1))}^2\leq C\,K^2\,R^{1+d+2\alpha}.
\end{align}

\emph{Step 4. Iteration.} We retain $y_1$ and define $\Theta\colon(0,1/2]\to\R_+$ by
\begin{align}\label{Eq:THETADEF}
\begin{split}
    \Theta(\varrho) &:= \int_{B_\varrho(y_1)}\big\vert\rmd u-(\rmd u)_{y_1,\varrho}\big\vert^2\d\Leb^{1+d}\\
    &\phantom{:}= \min\!\Big\{\!\int_{B_\varrho(y_1)}\vert\rmd u-k\vert^2\d\Leb^{1+d}: k\in\R^{1+d}\Big\},
    \end{split}
\end{align}
which is clearly nondecreasing. Let $\varrho,R\in (0,1/2]$ with $\varrho\leq R$ and let $v$ be from Step 3. By minimality of the mean, \eqref{Eq:Decay}, and \eqref{Eq:Comparison},
\begin{align*}
    \Theta(\varrho)&\leq \int_{B_\varrho(y_1)}\big\vert\rmd u- (\rmd v)_{y_1,\varrho}\big\vert^2\d\Leb^{1+d}\\
    &\leq 2\int_{B_\varrho(y_1)}\big\vert\rmd v-(\rmd v)_{y_1,\varrho}\big\vert^2\d\Leb^{1+d}+2\int_{B_\varrho(y_1)}\vert\rmd u-\rmd v\vert^2\d\Leb^{1+d}\\
    &\leq 2C\,\Big[\frac{\varrho}{R}\Big]^{1+d+2\gamma}\int_{B_R(y_1)}\big\vert\rmd v-(\rmd u)_{y_1,R}\big\vert^2\d\Leb^{1+d}+2C\,K^2\,R^{1+d+2\alpha}\\
    &\leq 4C\,\Big[\frac{\varrho}{R}\Big]^{1+d+2\gamma}\int_{B_R(y_1)}\big\vert\rmd u-(\rmd u)_{y_1,R}\big\vert^2\d\Leb^{1+d}\\
    &\qquad\qquad + 4C\int_{B_R(y_1)}\vert\rmd u-\rmd v\vert^2\d\Leb^{1+d}+ 2C^2\,K^2\,R^{1+d+2\alpha}\\
    &\leq 4C\,\Big[\frac{\varrho}{R}\Big]^{1+d+2\gamma}\,\Theta(R) + 6C^2\,K^2\,R^{1+d+2\alpha};
\end{align*}
here, we used $\varrho\leq R$ and $C\geq 1$ in the second last inequality.

Set $\beta:=\min\{\alpha,\gamma/2\}$. Since $R\leq 1$, we may replace $\smash{R^{1+d+2\alpha}}$ by $\smash{R^{1+d+2\beta}}$ in the last line. As $\beta<\gamma$, \cref{lem:iter} applies with $1/2$, $1+d+2\beta$, and $1+d+2\gamma$ in place of $R_0$, $a$, and $b$, respectively. This gives, for every $\varrho\in(0,1/2]$,
\begin{align*}
    \Theta(\varrho) &\leq C\,\varrho^{1+d+2\beta}\,\big[\Theta(1/2)+K^2\big]\leq C\,\varrho^{1+d+2\beta}\,\big[\big\Vert\rmd u-\zeta_0\big\Vert_{L^2(B_1(0))}^2+K^2\big];
\end{align*}
in the last step, we used $\zeta_0$ as a competitor in the definition \eqref{Eq:THETADEF} of $\Theta(1/2)$ and then the inclusion $B_{1/2}(y_1)\subset B_1(0)$.

\emph{Step 5. Conclusion.} Abbreviate $\smash{H:=\Vert\rmd u-\zeta_0\Vert_{L^2(B_1(0))}+K}$. As $\smash{y_1\in B_{1/2}(0)}$ was arbitrary in the previous Step 4, \cref{lem:campanato} with $z$ and $s$ replaced by $0$ and $1/2$, respectively, shows $\rmd u$ (and hence $\rmd u-k$ for every $k\in\R^{1+d}$) has a representative which is $\beta$-Hölder continuous on $\smash{\overline{B}_{1/2}(0)}$ with seminorm at most $CH$. In turn, for any given $y\in B_{1/2}(0)$ we estimate
\begin{align*}
    \vert\rmd u(y)-\zeta_0\vert&\leq\big\vert\rmd u(y)-(\rmd u)_{0,1/2}\big\vert+\big\vert(\rmd u)_{0,1/2}-\zeta_0\big\vert\\
    &\leq \frac{1}{\Leb^{1+d}[B_{1/2}(0)]}\int_{B_{1/2}(0)}\big\vert\rmd u(y) - \rmd u\big\vert\d\Leb^{1+d}\\
    &\qquad\qquad + \frac{1}{\Leb^{1+d}[B_{1/2}(0)]}\int_{B_{1/2}(0)}\vert \rmd u-\zeta_0\vert\d\Leb^{1+d} \\
    &\leq CH+\frac{1}{\Leb^{1+d}[B_{1/2}]}\,\Vert\rmd u-\zeta_0\Vert_{L^2(B_1(0))};
\end{align*}
the last inequality results from the previous Hölder estimate and $\diam B_{1/2}(0) = 1$ as well as the Cauchy--Schwarz inequality.
\end{proof}

\section{Regularity of $p$-harmonic temporal functions}\label{Sec:Reg}

In this section, we prove \cref{thm:reg}. Let $g$ be a Lorentzian metric on $M$ of class $\smash{C^{k,\alpha}_\loc}$, where $k\in\N_0$ and $\alpha\in(0,1)$. Finally, let $u\colon M\to\R$ be a $p$-harmonic temporal function on $M$.

The claim is local. Thus, let $x_0\in M$ and let an open ball $B\Subset M$ around $x_0$ and a compact convex set $\scrC$ be given by \cref{lem:K}. We then adopt the framework of \cref{Sec:Lemmas} for this choice of $B$ and $\scrC$. In particular, $c$, $\Lambda$, and $K$ denote the constants and $\smash{\overline{a}}$ denotes the truncation of $a$ from \cref{lem:trunc}, respectively. Let $\lambda_1$ and $\Lambda_1$ be the constants $\lambda$ and $\Lambda$ provided by \cref{lem:ell} for the compact set $\smash{\overline{B}\times\scrC\subset M\times \R^{1+d}}$. Since $\rmd u\in\scrC$ holds $\Leb^{1+d}$-a.e.~on $B$, by \eqref{Eq:OrigWeakSol} and item \ref{La:T1} from \cref{lem:trunc} the function $u$ belongs to $\smash{W^{1,\infty}(B)}$ and is a weak solution of
\begin{align}\label{Eq:Trunc}
    \partial_\alpha\big[\overline{a}^\alpha(\cdot,\rmd u)\big] = 0\quad\textnormal{in }B.
\end{align}
We will show $u$ is of class $\smash{C^{k+1,\alpha}_\loc}$ on $B$ by induction on $k$.

\begin{proof}[Proof of \cref{thm:reg}, $k=0$] Recall $u$ is Lipschitz, hence continuous.

To study its derivative, we proceed in two steps. We first show Hölder continuity of $\rmd u$ with some exponent $\beta\in(0,\alpha]$. Then we improve this exponent iteratively until we reach the Hölder exponent $\alpha$ of $g$.

\emph{Step 1. Hölder continuity of $\rmd u$ with some exponent.} Let $y\in B$ and $\varrho\in(0,1]$ with $B_\varrho(y)\subset B$. Define rescalings $\smash{\widetilde{u}\colon B_1(0)\to\R}$ and $\smash{\widetilde{a}\colon B_1(0)\times \R^{1+d}\to\R^{1+d}}$ of $u$ and $\overline{a}$ to $B_1(0)$, respectively, by
\begin{align*}
    \widetilde{u}(z) &:= \frac{u(y+\varrho z)}{\varrho},\\
    \widetilde{a}(z,\zeta) &:= \overline{a}(y+\varrho z,\zeta).
\end{align*}
Then $\smash{\rmd\widetilde u(z) = \rmd u(y+\varrho\,z)}$ for all $z\in B_1(0)$, $\smash{\widetilde u\in H^1(B_1(0))}$ is a weak solution of 
\begin{align*}
    \partial_\alpha\big[\widetilde{a}^\alpha(\cdot,\rmd\widetilde u)\big]=0\quad\textnormal{in }B_1(0),
\end{align*}
and by \cref{lem:trunc}, the map $\smash{\widetilde a}$ obeys \eqref{Eq:Struct} with the constants $c$, $\Lambda$, and $\varrho^\alpha\,K$ (and hence $K$). By \cref{lem:camp} and removing the scaling, $\rmd u$ has a representative which is $\beta$-Hölder continuous on $\smash{\overline{B}_{\varrho/2}(y)}$, where $\beta\in(0,\alpha]$ depends only on $d$, $c$, $\Lambda$, and $\alpha$. Since $y$ is arbitrary, this shows $\rmd u$ is of class $\smash{C^{0,\beta}_\loc}$ on $B$. In particular, $u\in C^1(B)$ and, as $\scrC$ is closed, the inclusion $\smash{\rmd u(x)\in \scrC}$ for $\smash{\Leb^{1+d}}$-a.e.~$x\in B$ from \cref{lem:K} upgrades to $\rmd u(x)\in\scrC$ for every $x\in B$. If $\beta = \alpha$, the proof is finished.

\emph{Step 2. Improvement of the exponent.} Otherwise, assume $\rmd u$ is of class $\smash{C^{0,\sigma}_\loc}$ on $B$ for some $\sigma\in(0,\alpha)$. We claim $\rmd u$ is of class $\smash{C^{0,\sigma'}_\loc}$ on $B$, where $\sigma' := \min\{\alpha,2\sigma\}$. Starting from the exponent $\beta$ from Step 1 and applying the claim finitely many times, we conclude $\rmd u$ is of class $\smash{C^{0,\alpha}_\loc}$ on $B$, as desired.

To show the claim, let $B'\Subset B$ be an open ball. Moreover, let $H\in\R_+$ denote the $\sigma$-Hölder seminorm of $\rmd u$ on $B'$. Let $x_1\in B'$ and $R\in(0,1]$ with $B_R(x_1)\subset B'$. Recall from Step 1 that $\rmd u(x_1)\in\scrC$. By \cref{lem:ell}, the symmetric matrix $\smash{\sfb^{\alpha\beta} := \partial_{\zeta_\beta}a^\alpha(x_1,\rmd u(x_1))}$ satisfies the hypotheses of \cref{lem:const} with the constants $\lambda_1$ and $\smash{(1+d)^{1/2}\,\Lambda_1}$ in place of $\lambda$ and $\Lambda$, respectively. For $x\in B_R(x_1)$, we set
\begin{align*}
    E(x) := a(x,\rmd u(x))-a(x_1,\rmd u(x_1))-\sfb\,\big[\rmd u(x)-\rmd u(x_1)\big].
\end{align*}
By \ref{La:T1} and \ref{La:T4} from \cref{lem:trunc}, we have 
\begin{align*}
\big\vert a(x,\rmd u(x))-a(x_1,\rmd u(x))\big\vert= \big\vert \overline{a}(x,\rmd u(x))-\overline{a}(x_1,\rmd u(x))\big\vert \leq K\,R^\alpha.  
\end{align*}
Moreover, since the segment from $\rmd u(x_1)$ to $\rmd u(x)$ is contained in the convex set $\scrC$, Taylor's theorem and $\sigma$-Hölder continuity of $\rmd u$ give
\begin{align*}
    &\big\vert a(x_1,\rmd u(x))-a(x_1,\rmd u(x_1))-\sfb\,\big[\rmd u(x)-\rmd u(x_1)\big]\big\vert\\
    &\qquad\qquad \leq C'\,\big\vert\rmd u(x)-\rmd u(x_1)\big\vert^2\\
    &\qquad\qquad \leq C'H^2\,R^{2\sigma},
\end{align*}
where $C'>0$ denotes the maximum of the norm of $\smash{\partial_\zeta^2a}$ on $\smash{\overline{B}\times\scrC}$. Combining these two estimates and using $R\leq 1$, we obtain
\begin{align}\label{Eq:ErrorBound}
    \Vert E\Vert_{C^0(B_R(x_1))}\leq\big[K+C'H^2\big]\,R^{\sigma'}.
\end{align}

Since $u$ is $p$-harmonic on $B_R(x_1)$ by \eqref{Eq:OrigWeakSol} and since constant vector fields are divergence-free, every $\smash{\varphi\in H^1_0(B_R(x_1))}$ satisfies
\begin{align}\label{Eq:Perturbed}
\begin{split}
    \int_{B_R(x_1)}\sfb^{\alpha\beta}\,\partial_\beta u\,\partial_\alpha\varphi\d\Leb^{1+d} &= \int_{B_R(x_1)} \big[\sfb^{\alpha\beta}\,\partial_\beta u-a^\alpha(\cdot,\rmd u)\big]\,\partial_\alpha\varphi\d\Leb^{1+d}\\
    &= -\int_{B_R(x_1)} E^\alpha\,\partial_\alpha\varphi\d\Leb^{1+d}.
    \end{split}
\end{align}
By the Lax--Milgram theorem, which is applicable thanks to \cref{lem:ell} and the Poincaré inequality, there exists a unique weak solution $\smash{v\in H^1(B_R(x_1))}$ of 
\begin{align*}
\partial_\alpha\big[\sfb^{\alpha\beta}\,\partial_\beta v\big]=0 \quad\textnormal{in }B_R(x_1)    
\end{align*}
with Dirichlet boundary conditions $\smash{v-u\in H^1_0(B_R(x_1))}$. Testing this equation and \eqref{Eq:Perturbed} with $u-v$ and employing \eqref{Eq:ErrorBound}, we get
\begin{align*}
    \lambda_1\,\big\Vert\rmd u-\rmd v\big\Vert_{L^2(B_R(x_1))}^2 &\leq -\int_{B_R(x_1)}E^\alpha \,\partial_\alpha(u-v)\d\Leb^{1+d}\\
    &\leq\big[K+C'H^2\big]\,R^{\sigma'}\int_{B_R(x_1)}\vert\rmd u-\rmd v\vert\d\Leb^{1+d}.
\end{align*}
By the Cauchy--Schwarz inequality, this implies
\begin{align*}
    \big\Vert\rmd u-\rmd v\big\Vert_{L^2(B_R(x_1))}^2\leq C\,\big[K+C'H^2\big]^2\,R^{1+d+2\sigma'},
\end{align*}
where the constant $C>0$ depends only on $d$ and $\lambda_1$. 

Given any $\varrho >0$ such that $\smash{B_\varrho(x_1)\subset B'}$,  as in \eqref{Eq:THETADEF} we set
\begin{align*}
    \Theta(\varrho) := \int_{B_\varrho(x_1)}\big\vert\rmd u-(\rmd u)_{x_1,\varrho}\big\vert^2\d\Leb^{1+d}.
\end{align*}
Arguing as in Step 4 of the proof of \cref{lem:camp}, with \cref{lem:const} in place of \eqref{Eq:Decay}, cf.~also the proof of Giaquinta--Martinazzi \cite[Thm.~5.19]{giaquinta-martinazzi2012}, we infer there exists a constant $C>0$ depending only on $d$, $\lambda_1$, and $\Lambda_1$ such that for every $\varrho\in (0,R]$,
\begin{align*}
    \Theta(\varrho)\leq C\,\Big[\frac{\varrho}{R}\Big]^{3+d}\,\Theta(R)+C\,\big[K+C'H^2\big]^2\,R^{1+d+2\sigma'}.
\end{align*}

Now let $B''\Subset B'$ be a concentric open ball. We let $R_0\in(0,1]$ be at most half the difference of the radii of $B'$ and $B''$. Then $B_{R_0}(x_1)\subset B'$ for every $x_1$ whose distance to $B''$ is less than $R_0$. Since $\sigma'\leq\alpha<1$, \cref{lem:iter} applies with the exponents $1+d+2\sigma'$ and $3+d$ in place of $a$ and $b$, respectively, and shows that for every such $x_1$ and every $\varrho\in(0,R_0]$ we have
\begin{align*}
    \Theta(\varrho)\leq C\,\varrho^{1+d+2\sigma'}\Big[R_0^{-1-d-2\sigma'}\int_{B'}\vert\rmd u\vert^2\d\Leb^{1+d}+\big[K+C'H^2\big]^2\Big].
\end{align*}
By \cref{lem:campanato}, applied with $R_0$ in place of $s$ and an arbitrary center $z\in B''$, we see $\rmd u$ is $\sigma'$-Hölder continuous on a neighborhood of every point of $B''$. Since $B''$ and $ B'$ were arbitrary, the claim follows.
\end{proof}

\begin{proof}[Proof of \cref{thm:reg}, $k\in\N$] We argue by induction on $k\in\N_0$. The previous proof establishes the claim for $k=0$. Assume the claim of \cref{thm:reg} holds with $k-1$ in place of $k$, where $k\in\N$. As $g$ is of class $\smash{C^{k,\alpha}_\loc}$, hence of class $\smash{C^{k-1,\alpha}_\loc}$, the function $u$ is of class $\smash{C^{k,\alpha}_\loc}$ on $B$ by the induction hypothesis. In particular, as above we infer $u\in C^1(B)$ and $\rmd u(x)\in\scrC$ for every $x\in B$. By \eqref{Eq:DefA}, the map $a$ is of class $\smash{C^{k,\alpha}}$ in its first argument and smooth in its second argument on a neighborhood of $\smash{\overline{B}\times\scrC}$, together with all derivatives in the second argument.

\emph{Step 1. Second derivatives.} We first claim $\smash{u\in H^2_\loc(B)}$. Let $B''\Subset B'\Subset B$ be concentric open balls,  let $\smash{e_\kappa\in\R^{1+d}}$ be a coordinate direction, and let $h\in\R\setminus\{0\}$ such that $\vert h\vert$ is smaller than the distance of $B'$ to $\partial B$. Define $\smash{z_h\in H^1(B')}$ by
\begin{align*}
    z_h := \frac{u(\cdot +he_\kappa)-u}{h}.
\end{align*}
For $\varphi\in H_0^1(B')$, we test the $p$-harmonicity \eqref{Eq:OrigWeakSol} of $u$ with $\varphi(\cdot-he_\kappa)$ and $\varphi$ to get
\begin{align}\label{Eq:PDEEEse}
    \int_{B'}\big[a^\alpha(\cdot+he_\kappa,\rmd u(\cdot +he_\kappa))-a^\alpha(\cdot,\rmd u)\big]\,\partial_\alpha\varphi\d\Leb^{1+d} = 0.
\end{align}
As $\rmd u(\cdot+he_\kappa)-\rmd u = h\,\rmd z_h$, the term in brackets equals $\smash{h\,[\sfb_h^{\alpha\beta}\,\partial_\beta z_h+F_h^\alpha]}$, where
\begin{align*}
    \sfb_h^{\alpha\beta}(x) &:= \int_{[0,1]}\partial_{\zeta_\beta}a^\alpha(x+he_\kappa,\rmd u(x)+th\,\rmd z_h(x))\d t,\\
    F_h^\alpha(x) &:= \frac{a^\alpha(x+he_\kappa,\rmd u(x))-a^\alpha(x,\rmd u(x))}{h}.
\end{align*}

By convexity of $\scrC$ and \cref{lem:ell}, every $\smash{\xi\in\R^{1+d}}$ satisfies 
\begin{align*}
\sfb_h^{\alpha\beta}\,\xi_\alpha\xi_\beta &\geq\lambda_1\,\vert\xi\vert^2,\\
\vert\sfb_h\vert&\leq(1+d)^{1/2}\,\Lambda_1.    
\end{align*}
Since $k\geq 1$, $\vert F_h\vert$ is bounded by the maximum $C'$ of the norm of $\smash{\partial_x a}$ on $\smash{\overline{B}\times\scrC}$. Let $\smash{\eta\in C^\infty_\comp(B')}$ be a cutoff function with $\eta=1$ on $B''$. Replacing $\varphi$ by $\smash{\eta^2z_h}$ in \eqref{Eq:PDEEEse}, which belongs to $\smash{H_0^1(B')}$, and dividing by $h$ yields
\begin{align*}
    0 &= \int_{B'} \big[\sfb_h^{\alpha\beta}\,\partial_\beta z_h + F_h^\alpha\big]\,\partial_\alpha(\eta^2z_h)\d\Leb^{1+d}\\
    &= 2\int_{B'}\eta\, z_h\,\sfb_h^{\alpha\beta}\,\partial_\beta z_h\,\partial_\alpha\eta\d\Leb^{1+d} + \int_{B'}\eta^2\,\sfb_h^{\alpha\beta}\,\partial_\beta z_h\,\partial_\alpha z_h\d\Leb^{1+d}\\
    &\qquad\qquad +\int_{B'} F_h^\alpha\,\big[\eta^2\,\partial_\alpha z_h + 2\eta\,z_h\,\partial_\alpha\eta\big]\d\Leb^{1+d}.
\end{align*}
Consequently,
\begin{align*}
    \lambda_1\int_{B'}\eta^2\,\vert\rmd z_h\vert^2\d\Leb^{1+d}&\leq\int_{B'}\eta^2\,\sfb_h^{\alpha\beta}\,\partial_\beta z_h\,\partial_\alpha z_h\d\Leb^{1+d}\\
    &= -\int_{B'}2\,\eta\,z_h\,\sfb_h^{\alpha\beta}\,\partial_\beta z_h\,\partial_\alpha\eta\d\Leb^{1+d}\\
    &\qquad\qquad-\int_{B'}F_h^\alpha\,\big[\eta^2\,\partial_\alpha z_h+2\,\eta\,z_h\,\partial_\alpha\eta\big]\d\Leb^{1+d}.
\end{align*}
We estimate the three terms on the right-hand side separately. Let $L>0$ be a Lipschitz constant of $u$ on $B$, whence $\vert z_h\vert\leq L$ on $B'$. Together with the inequalities $\smash{\vert\sfb_h\vert\leq(1+d)^{1/2}\Lambda_1}$ and $\vert F_h\vert\leq C'$ as well as Young's inequality, we easily infer
\begin{align*}
    -\int_{B'}2\,\eta\,z_h\,\sfb_h^{\alpha\beta}\,\partial_\beta z_h\,\partial_\alpha\eta\d\Leb^{1+d} &\leq\frac{\lambda_1}{4}\int_{B'}\eta^2\,\vert\rmd z_h\vert^2\d\Leb^{1+d}\\
    &\qquad\qquad +\frac{4(1+d)\,\Lambda_1^2\,L^2}{\lambda_1}\int_{B'}\vert\rmd\eta\vert^2\d\Leb^{1+d},\\
    -\int_{B'}F_h^\alpha\,\eta^2\,\partial_\alpha z_h\d\Leb^{1+d} &\leq\frac{\lambda_1}{4}\int_{B'}\eta^2\,\vert\rmd z_h\vert^2\d\Leb^{1+d}+\frac{C'^2}{\lambda_1}\int_{B'}\eta^2\d\Leb^{1+d},\\
    -2\int_{B'}F_h^\alpha\,\eta\,z_h\,\partial_\alpha\eta\d\Leb^{1+d} &\leq 2C'L\int_{B'}\eta\,\vert\rmd\eta\vert\d\Leb^{1+d}.
\end{align*}
Combining these estimates and rearranging terms, 
\begin{align*}
    \frac{\lambda_1}{2} \big\Vert\rmd z_h\big\Vert_{L^2(B'')}^2 \leq \frac{\lambda_1}{2}\int_{B'}\eta^2\,\vert\rmd z_h\vert^2\d\Leb^{1+d} \leq C,
\end{align*}
where $C>0$ is a constant independent of $h$. Thus, the difference quotients of $\rmd u$ are bounded in $\smash{L^2(B'')}$ uniformly in $h$. By Gilbarg--Trudinger \cite[Lem.~7.24]{gilbarg-trudinger2001}, this yields $\smash{u\in H^2(B'')}$, which proves the claim since $B''$ was arbitrary.

\emph{Step 2. Differentiated equation.} By Step 1 and the chain rule, the vector field $a(\cdot,\rmd u)$ belongs to $\smash{H^1_\loc(B)}$ and satisfies
\begin{align*}
    \partial_\kappa\big[a^\alpha(\cdot,\rmd u)\big] = \sfb^{\alpha\beta}\,\partial_\beta\partial_\kappa u+F^\alpha\quad\Leb^{1+d}\mres B\textnormal{-a.e.},
\end{align*}
where we abbreviate
\begin{align*}
    \sfb^{\alpha\beta} &:= \partial_{\zeta_\beta}a^\alpha(\cdot,\rmd u),\\
    F^\alpha &:= \partial_\kappa a^\alpha(\cdot,\rmd u).
\end{align*}
Testing the PDE \eqref{Eq:OrigWeakSol} certifying $p$-harmonicity of $u$ with $\partial_\kappa\varphi$, where $\smash{\varphi\in C_\comp^\infty(B)}$, and integrating by parts, we see 
\begin{align*}
    0 &= \int_B a^\alpha(\cdot,\rmd u)\,\partial_\alpha\partial_\kappa\varphi\d\Leb^{1+d}\\
    &= -\int_B \partial_\kappa\big[a^\alpha(\cdot,\rmd u)\big]\,\partial_\alpha\varphi\d\Leb^{1+d}\\
    &= -\int_B \big[\sfb^{\alpha\beta}\,\partial_\beta\partial_\kappa u + F^\alpha\big]\,\partial_\alpha\varphi\d\Leb^{1+d}.
\end{align*}
Since $\varphi$ is arbitrary, this means the function $\smash{w\in H_\loc^1(B)}$ given by $\smash{w:=\partial_\kappa u}$ is a weak solution of the linear equation
\begin{align*}
    \partial_\alpha\big[\sfb^{\alpha\beta}\,\partial_\beta w\big] = -\partial_\alpha F^\alpha\quad\textnormal{in }B.
\end{align*}
Since $\rmd u$ is of class $\smash{C^{k-1,\alpha}_\loc}$ on $B$ with values in $\scrC$, both $\sfb$ and $F$ are of class $\smash{C^{k-1,\alpha}_\loc}$ on $B$ and $\sfb$ is uniformly elliptic thanks to \cref{lem:ell}.

\emph{Step 3. Conclusion.} By Schauder theory for linear elliptic PDEs in divergence form, cf.~Giaquinta--Martinazzi \cite[Thms.~5.19, 5.20]{giaquinta-martinazzi2012} or Gilbarg--Trudinger \cite[Thm.~8.32]{gilbarg-trudinger2001}, the function $w$ from Step 2 is of class $\smash{C^{k,\alpha}_\loc}$ on $B$. Since $\kappa$ was arbitrary, this shows $u$ is of class $\smash{C^{k+1,\alpha}_\loc}$ on $B$. This completes the proof.
\end{proof}

For later use, we save a variation of Step 1 for locally Lipschitz metrics.

\begin{lemma}[Second derivatives for $p$-harmonic temporal functions]\label{lem:H2} Assume $g$ is locally Lipschitz. Then every $p$-harmonic temporal function on $M$ lies in $\smash{H^2_\loc}$.
\end{lemma}

\begin{proof} As $g$ is of class $\smash{C^{0,\sigma}_\loc}$ for every $\sigma\in(0,1)$, \cref{thm:reg} gives $u\in C^1(M)$. Let $x_0\in M$ and let $B$, $\scrC$, $\lambda_1$, and $\Lambda_1$ be as above; as in Step 1 of the proof of \cref{thm:reg} for $k=0$, we obtain $\rmd u(x)\in\scrC$ for every $x\in B$. Then Step 1 of the proof of \cref{thm:reg} for $k\in\N$ carries over almost verbatim and yields $\smash{u\in H^2_\loc(B)}$. Indeed, retaining the notation therein, the only place where that step uses the regularity of $g$ is the bound on $F_h$. Under the present hypothesis, the argument for item \ref{La:T4} from \cref{lem:trunc} with $\alpha$ replaced by $1$ shows $a(\cdot,\zeta)$ is Lipschitz on $\smash{\overline{B}}$ uniformly in $\zeta\in\scrC$. Hence, $\vert F_h\vert$ is bounded by the corresponding Lipschitz constant in place of the maximum of the norm of $\partial_xa$ on $\smash{\overline{B}\times\scrC}$. As $x_0$ is arbitrary, the claim follows.
\end{proof}

\section{Existence of $p$-harmonic coordinates}\label{Sec:Exist}

In this section, we prove \cref{thm:coord}. Let $g$ be a Lorentzian metric on $M$ of class $\smash{C^{k,\alpha}_\loc}$, where $k\in\N_0$ and $\alpha\in(0,1)$. 

The core of the argument is the following local existence result for $p$-harmonic temporal functions with prescribed differential at a point, up to a small error.

\begin{proposition}[Local existence with tilt estimate]\label{prop:exist} Let $x_0\in M$. Let $\smash{\omega\in\R^{1+d}}$ be a future-directed timelike covector at $x_0$. Then there exist $r_0>0$ and a constant $C>0$ such that for every $r\in(0,r_0)$, there exists a $p$-harmonic temporal function $\smash{\tau\in C^1(B_r(x_0))}$ such that
\begin{align*}
    \Vert\rmd\tau-\omega\Vert_{C^0(B_r(x_0))}\leq C\,r^\alpha.
\end{align*}
\end{proposition}

\begin{proof} \emph{Step 1. Truncation.} The set of all $\smash{(x,\zeta)\in M\times\R^{1+d}}$ such that $\zeta$ is a future-directed timelike covector at $x$ is open by continuity of $g$ and the time-orientation $X$. Therefore, there exist $\delta>0$ and an open ball $B\Subset M$ around $x_0$ such that for every $\smash{x\in\overline{B}}$, every element of the closed ball $\smash{\scrC := \overline{B}_\delta(\omega)}$ is a future-directed timelike covector at $x$. We adopt the framework of \cref{Sec:Lemmas} for this choice of $B$ and $\scrC$. In particular, $c$, $\Lambda$, and $K$ denote the constants and $\smash{\overline{a}}$ denotes the truncation of $a$ from \cref{lem:trunc}. Finally, define $\smash{\ell\colon \R^{1+d}\to\R}$ by $\ell(x) := \omega(x-x_0)$. Note that if $\omega$ were the differential of a temporal function at $x_0$, then $\ell$ would be the linearization of that function at $x_0$ up to an additive constant.

\emph{Step 2. Dirichlet problem.} Let $r>0$ with $\smash{B_{2r}(x_0)\subset B}$. By items \ref{La:T2} and \ref{La:T3} from \cref{lem:trunc} and the Poincaré inequality, the assignment
\begin{align*}
    \langle T\chi,\varphi\rangle := \int_{B_{2r}(x_0)}\overline{a}^\alpha(\cdot,\omega+\rmd\chi)\,\partial_\alpha\varphi\d\Leb^{1+d}
\end{align*}
for every $\smash{\varphi\in H_0^1(B_{2r}(x_0))}$ does define a Lipschitz and strongly monotone operator $T\colon H^1_0(B_{2r}(x_0))\to H^{-1}(B_{2r}(x_0))$. As in Step 3 in the proof of \cref{lem:camp}, Zarantonello's theorem \cite[Thm.~25.B]{zeidler1990b} implies there is a unique $\smash{\chi\in H^1_0(B_{2r}(x_0))}$ with $T\chi=0$. This shows the function $\smash{\tau\in H^1(B_{2r}(x_0))}$ given by $\tau := \ell+\chi$ satisfies Dirichlet boundary conditions $\smash{\tau-\ell\in H^1_0(B_{2r}(x_0))}$ and is a weak solution of
\begin{align}\label{Eq:TruncTau}
    \partial_\alpha\big[\overline{a}^\alpha(\cdot,\rmd\tau)\big]=0\quad\textnormal{in }B_{2r}(x_0).
\end{align}

\emph{Step 3. Rescaling and closeness in $L^2$.} We define the rescalings $\smash{u\colon B_1(0)\to\R}$ and $\smash{\overline{a}_r\colon B_1(0)\times \R^{1+d}\to \R^{1+d}}$ of $\tau$ and $\overline{a}$ to $B_1(0)$, respectively, by
\begin{align*}
    u(y) &:= \frac{\tau(x_0+2r\,y)}{2r},\\
    \overline{a}_r(y,\zeta) &:= \overline{a}(x_0+2r\,y,\zeta).
\end{align*}
Then we have $\rmd u(y) = \rmd\tau(x_0+2r\,y)$ for every $y\in B_1(0)$. Identifying $\omega$ with a linear function, we have $\smash{u-\omega\in H^1_0(B_1(0))}$; by \eqref{Eq:TruncTau}, $u$ is a weak solution of 
\begin{align*}
\partial_\alpha\big[\overline{a}_r^\alpha(\cdot,\rmd u)\big]=0\quad\textnormal{in }B_1(0).    
\end{align*}
By \cref{lem:trunc}, the map $\overline{a}_r$ satisfies \eqref{Eq:Struct} with the constants $c$, $\Lambda$, and $(2r)^\alpha\,K$. Since $\smash{u-\omega\in H^1_0(B_1(0))}$, testing the previous equation with $u-\omega$ gives
\begin{align*}
    \int_{B_1(0)}\overline{a}_r^\alpha(\cdot,\rmd u)\,\partial_\alpha(u-\omega)\d\Leb^{1+d} = 0.
\end{align*}
Since $\smash{\overline{a}_r(0,\omega)}$ is a constant vector,
\begin{align*}
    \int_{B_1(0)}\overline{a}_r^\alpha(0,\omega)\,\partial_\alpha(u-\omega)\d\Leb^{1+d} = 0.
\end{align*}
Together with \eqref{Eq:Struct}, these two identities yield
\begin{align*}
    c\,\big\Vert\rmd u-\omega\big\Vert_{L^2(B_1(0))}^2  &\leq\int_{B_1(0)}\big[\overline{a}_r^\alpha(\cdot,\rmd u)-\overline{a}_r^\alpha(\cdot,\omega)\big]\,\partial_\alpha(u-\omega)\d\Leb^{1+d}\\
    &= -\int_{B_1(0)}\overline{a}_r^\alpha(\cdot,\omega)\,\partial_\alpha(u-\omega)\d\Leb^{1+d}\\
    &= -\int_{B_1(0)}\big[\overline{a}_r^\alpha(\cdot,\omega)-\overline{a}_r^\alpha(0,\omega)\big]\,\partial_\alpha(u-\omega)\d\Leb^{1+d}\\
    &\leq (2r)^\alpha\,K\,\int_{B_1(0)} \vert\rmd u-\omega\vert\d\Leb^{1+d}.
\end{align*}
By the Cauchy--Schwarz inequality, this shows
\begin{align*}
    \Vert\rmd u-\omega\Vert_{L^2(B_1(0))}\leq C\,r^\alpha,
\end{align*}
where the constant $C>0$ depends only on $d$, $c$, $K$, and $\alpha$.

\emph{Step 4. Uniform closeness and removal of the truncation.} By \cref{lem:camp} with $\zeta_0$ replaced by $\omega$ and Step 3, $\rmd u$ is continuous on $\smash{\overline{B}_{1/2}(0)}$ with
\begin{align*}
    \Vert\rmd u-\omega\Vert_{C^0(B_{1/2}(0))}\leq C\,\big[\Vert\rmd u-\omega\Vert_{L^2(B_1(0))}+(2r)^\alpha\,K\big]\leq C'\,r^\alpha,
\end{align*}
where $C'>0$ is a constant depending only on $d$, $c$, $\Lambda$, $K$, and $\alpha$. Let $r_0>0$ be such that $\smash{B_{2r_0}(x_0)\subset B}$ and $\smash{C'\,r_0^\alpha\leq\delta}$ and let $r\in(0,r_0)$. Then $\tau\in C^1(B_r(x_0))$ and $\rmd\tau(x)\in\scrC$ for every $x\in B_r(x_0)$ by the above $C^0$-estimate. In particular, $\tau$ is Lipschitz continuous on $B_r(x_0)$. Since $\smash{\overline{B}\times\scrC}$ is compact and since $\zeta$ is a future-directed timelike covector at $x$ for every $\smash{(x,\zeta)\in\overline{B}\times\scrC}$, the continuous and positive quantities $\smash{-g^{-1}(x)(\zeta,\zeta)}$ and $\zeta(X(x))$ are both bounded from below by a positive constant on $\smash{\overline{B}\times\scrC}$. Hence, $\tau$ is a temporal function on $B_r(x_0)$ in the sense of \cref{Def:Temporal}. By \ref{La:T1} from \cref{lem:trunc} and \eqref{Eq:TruncTau}, we have $\smash{\Box_p\tau=0}$ in $\smash{\scrD'(B_r(x_0))}$. The claimed estimate follows by removing the rescaling.
\end{proof}

\begin{proof}[Proof of \cref{thm:coord}] As the set of all future-directed timelike covectors at $x_0$ is a nonempty open subset of $\smash{\R^{1+d}}$, it contains $1+d$ linearly independent covectors $\smash{\omega^0,\dots,\omega^d}$. Let $\smash{W\in\R^{(1+d)\times(1+d)}}$ be the invertible matrix with rows $\smash{\omega^0,\dots,\omega^d}$. Observe that $\vert W\,v\vert\geq \vert W^{-1}\vert^{-1}\,\vert v\vert$ for every $\smash{v\in\R^{1+d}}$. Let $r_0$ and $C$ be the minimum and the maximum, respectively, of the corresponding constants from the $(1+d)$-fold application of \cref{prop:exist} to $\smash{\omega^0,\dots,\omega^d}$. Let $r\in (0,r_0/2)$ satisfy $\smash{(1+d)^{1/2}\,C\,(2r)^\alpha\leq \vert W^{-1}\vert^{-1}/2}$. Let $\smash{\tau^0,\dots,\tau^d}$ designate the $p$-harmonic temporal functions on $B_{2r}(x_0)$ given by \cref{prop:exist}. Finally, define $\Phi\colon B_{2r}(x_0)\to\R^{1+d}$ by $\smash{\Phi:=(\tau^0,\dots,\tau^d)}$. By \cref{prop:exist}, $\Phi$ is of class $C^1$ and obeys the estimate $\vert\rmd\Phi(x)-W\vert\leq \vert W^{-1}\vert^{-1}/2$ for every $x\in B_{2r}(x_0)$ by our choice of $r$. For every $x,y\in B_{2r}(x_0)$, our preceding observations imply
\begin{align*}
    \big\vert\Phi(y)-\Phi(x)\big\vert &= \Big\vert\!\int_{[0,1]}\rmd\Phi(x+t\,(y-x))\,(y-x)\d t\Big\vert\\
    &= \Big\vert W(y-x)+ \int_{[0,1]}\big[\rmd\Phi(x+t\,(y-x))-W\big]\,(y-x)\d t\Big\vert\\
    &\geq \big\vert W(y-x)\big\vert - \int_{[0,1]}\big\vert\rmd\Phi(x+t\,(y-x))-W\big\vert\,\vert y-x\vert\d t\\
    &\geq \vert W^{-1}\vert^{-1}\,\vert y-x\vert - \frac{\vert W^{-1}\vert^{-1}}{2}\,\vert y-x\vert\\
    &= \frac{\vert W^{-1}\vert^{-1}}{2}\,\vert y-x\vert.
\end{align*}
Hence, $\Phi$ is injective on $B_{2r}(x_0)$. Moreover, dividing the previous estimate with $x+tv$ in place of $y$, where $\smash{v\in\R^{1+d}}$, by $t>0$ and letting $t\to 0$ yields $\vert\rmd\Phi(x)\,v\vert\geq\vert W^{-1}\vert^{-1}\,\vert v\vert/2$ for every $x\in B_{2r}(x_0)$; in particular, $\rmd\Phi(x)$ is invertible. By the inverse function theorem, $\Phi$ is a local $C^1$-diffeomorphism. Since $\Phi$ is injective, it is a $C^1$-diffeomorphism from $B_{2r}(x_0)$ onto its image, which is open. In fact, by \cref{thm:reg}, $\Phi$ is of class $\smash{C^{k+1,\alpha}_\loc}$ on $B_{2r}(x_0)$, and so is $\Phi^{-1}$ on $\Phi(B_{2r}(x_0))$. Since $\smash{\overline{B}_r(x_0)}$ and its image under $\Phi$ are compact, the restriction of $\Phi$ to $B_r(x_0)$ is a $\smash{C^{k+1,\alpha}}$-diffeomorphism onto its image and the restrictions of $\smash{\tau^0,\dots,\tau^d}$ to $B_r(x_0)$ are $p$-harmonic temporal functions. 

The coefficients of $\Phi_*g$ are given by
\begin{align*}
    (\Phi_*g)_{\mu\nu} = g_{\alpha\beta}\circ\Phi^{-1}\,\partial_\mu(\Phi^{-1})^\alpha\,\partial_\nu(\Phi^{-1})^\beta;
\end{align*}
the right-hand side is of class $\smash{C^{k,\alpha}}$.

Finally, assume $g$ is smooth. Then $g$ is of class $\smash{C^{k,1/2}_\loc}$ for every $k\in\N_0$. Since the construction of $\Phi$ and of its domain does not depend on $k$, the above argument shows that $\Phi$ and $\smash{\Phi^{-1}}$ are of class $\smash{C^{k+1,1/2}}$ for every $k\in\N_0$, hence smooth.
\end{proof}

\section{Applications}\label{Sub:RegSptIso}

\subsection{Invariances} Let $g$ and $g'$ be continuous Lorentzian metrics on open subsets $\smash{M,M'\subset\R^{1+d}}$, time-oriented by continuous vector fields $X$ and $X'$ on $M$ and $M'$, respectively. All notions  relative to $g'$ are marked with a prime. Let $\psi\colon M\to M'$ be locally Lipschitz. By Rademacher's theorem, $\psi$ is differentiable $\smash{\Leb^{1+d}}$-a.e.

\begin{definition}[Isometry]\label{Def:Isometry} Let $x\in M$ be a point where $\psi$ is differentiable. We say $\psi$ is an \textbf{isometry} at $x$ if    
\begin{align}\label{Eq:AeIso}
        g'_{\gamma\delta}\circ\psi(x)\,\rmd\psi^\gamma(x)\otimes\rmd \psi^\delta(x) = g(x);
    \end{align}
    we say $\psi$ is an isometry if it is an isometry $\vol$-a.e.
\end{definition}

\begin{definition}[Time-orientation-preservation] Let $x\in M$ be a point where $\psi$ is differentiable. We say $\psi$ is \textbf{time-orientation-preserving} at $x$ if its differential $\rmd\psi(x)$ maps future-directed $g(x)$-timelike vectors at $x$ to future-directed $g'\circ\psi(x)$-timelike vectors at $\psi(x)$; we say that $\psi$ is time-orientation-preserving if it is time-orientation-preserving $\vol$-a.e.
\end{definition}

\begin{lemma}[Invariance of temporality]\label{lem:inv} Let $\psi\colon M\to M'$ be a locally bi-Lipschitz map which is a time-orientation-preserving isometry. If $u'\colon M'\to\R$ is temporal relative to $g'$, then $u'\circ\psi$ is temporal relative to $g$.
\end{lemma}

\begin{proof} We abbreviate $u:=u'\circ\psi$, which is locally Lipschitz. 

As $\psi$ and $\smash{\psi^{-1}}$ are locally Lipschitz, both maps send negligible sets to negligible sets with respect to $\smash{\Leb^{1+d}}$. By \cref{Def:Temporal}, applied to a countable exhaustion of $M'$ by compact sets, $\rmd u'$ is a future-directed $g'$-timelike covector $\vol'$-a.e. Thus, there exists an $\smash{\Leb^{1+d}}$-negligible set $Z\subset M$ such that for every $x\in M\setminus Z$, the map $\psi$ is differentiable at $x$, $\psi$ is a time-orientation-preserving isometry at $x$, the function $u'$ is differentiable at $\psi(x)$, and $\rmd u'(\psi(x))$ is a future-directed $g'\circ\psi(x)$-timelike covector at $\psi(x)$. Let $K\subset M$ be  a compact set. Let $c'$ denote the constant from \cref{Def:Temporal} for $u'$ for the compact set $\psi(K)\subset M'$. Let $x\in K$ satisfy the previous properties. Moreover, we may and will assume $\psi(x)$ satisfies \eqref{Eq:LowBo}. 

By the chain rule, $u$ is differentiable at $x$ with $\rmd u(x) = \rmd u'\circ\psi(x)\,\rmd\psi(x)$. Taking the determinant of \eqref{Eq:AeIso}, we infer the matrix $\rmd\psi(x)$ is invertible. Inverting \eqref{Eq:AeIso}, we see the map $\smash{T\colon T_{\psi(x)}^*M'\to T_x^*M}$ with $\smash{(T\zeta')_\mu := \zeta'_\gamma\, {\rmd \psi(x)^\gamma}_\mu}$ is a linear isometry between $\smash{g'^{-1}\circ\psi(x)}$ and $\smash{g^{-1}(x)}$; in particular,
\begin{align}\label{Eq:NormInv}
    -g^{-1}(x)(\rmd u(x),\rmd u(x)) = -g'^{-1}\circ\psi(x)(\rmd u'\circ\psi(x),\rmd u'\circ\psi(x))\geq c'.
\end{align}

Since $\psi$ is a time-orientation-preserving isometry at $x$, the vector $\rmd\psi(x)\,X(x)$ is future-directed $g'\circ\psi(x)$-timelike and its length with respect to $g'\circ\psi(x)$ equals $\smash{\vert X(x)\vert_{g(x)}}$. By the chain rule and the reverse Cauchy--Schwarz inequality,
\begin{align*}
    \rmd u(x)(X(x)) = \rmd u'\circ\psi(x)(\rmd\psi(x)\,X(x))\geq\big\vert\rmd u'\big\vert_{g'}\circ \psi(x)\,\big\vert X(x)\big\vert_{g(x)}.
\end{align*}
The continuous and positive quantity $\smash{\vert X\vert_g}$ is bounded away from zero on $K$, say by $a>0$. Together with \eqref{Eq:NormInv}, this shows $\smash{\rmd u(x)(X(x)) \geq c'^{1/2}\,a}$, as desired.
\end{proof}

\begin{lemma}[Invariance of $p$-harmonicity]\label{Le:Invarp} Let $\psi\colon M\to M'$ be a locally bi-Lipschitz map that is a time-orientation-preserving isometry. If $u'\colon M'\to\R$ is a $p$-harmonic temporal function for $g'$, then $u'\circ\psi$ is a $p$-harmonic temporal function for $g$.
\end{lemma}

\begin{proof} Abbreviate $u:= u'\circ\psi$. By \cref{lem:inv}, $u$ is temporal. We also retain the notation from the previous proof.

It remains to show $\smash{\Box_pu=0}$ in $\scrD'(M)$. Let $\smash{\varphi\in C_\comp^\infty(M)}$. By a smooth partition of unity subordinate to a finite cover of $\supp\varphi$ and linearity of $\langle\Box_pu,\cdot\rangle$, we may and will assume $\supp\varphi\subset U$ for an open set $U\subset M$ such that $\psi$ is bi-Lipschitz on $U$; by invariance of domain, $\psi(U)$ is open and $\psi$ restricts to a homeomorphism from $U$ onto $\psi(U)$. The function $\varphi'\colon M'\to\R$ given by $\smash{\varphi' := \varphi\circ\psi^{-1}}$ on $\psi(U)$ and $\varphi':=0$ otherwise is Lipschitz with compact support in $M'$. Enlarging $Z$, we may and will assume $\varphi'$ is differentiable at $\psi(x)$ for every $x\in U\setminus Z$, entailing the chain rule $\rmd\varphi(x) = \rmd\varphi'\circ \psi(x)\,\rmd\psi(x)$. By isometry of $\psi$ at $x$  and \eqref{Eq:NormInv},
\begin{align*}
    \big\vert\rmd u\big\vert_g^{p-2}\,g^{-1}(\rmd u,\rmd\varphi) = \big[\big\vert\rmd u'\big\vert_{g'}^{p-2}\,g'^{-1}(\rmd u',\rmd\varphi')\big]\circ\psi(x).
\end{align*}
Moreover, taking the determinant of \eqref{Eq:AeIso} implies 
\begin{align*}
\big\vert\!\det\rmd\psi(x)\big\vert = \frac{\big\vert\!\det g(x)\big\vert^{1/2}}{\big\vert\!\det g'\circ\psi(x)\big\vert^{1/2}}.
\end{align*}
Thus, by the change of variables formula for locally Lipschitz homeomorphisms, cf. e.g. Evans--Gariepy \cite[p.~99]{evans-gariepy1992}, we infer
\begin{align*}
    -\big\langle\Box_p u,\varphi\big\rangle&=\int_U\big\vert\rmd u\big\vert_g^{p-2}\,g^{-1}(\rmd u,\rmd\varphi)\,\vert\!\det g\vert^{1/2}\d\Leb^{1+d}\\
    &= \int_U\big[\big\vert\rmd u'\big\vert_{g'}^{p-2}\,g'^{-1}(\rmd u',\rmd\varphi')\,\vert\!\det g'\vert^{1/2}\big]\circ\psi\,\vert\!\det\rmd\psi\vert\d\Leb^{1+d}\\
    &= \int_{\psi(U)}\big\vert\rmd u'\big\vert_{g'}^{p-2}\,g'^{-1}(\rmd u',\rmd\varphi')\,\vert\!\det g'\vert^{1/2}\d\Leb^{1+d}\\
    &= 0,
\end{align*}
where the last identity holds since $u'$ is $p$-harmonic and because  \eqref{Eq:Weak} extends to Lipschitz test functions with compact support.
\end{proof}

As a first consequence, we justify the gauge condition \eqref{Eq:Gauge}.

\begin{corollary}[Gauge condition]\label{cor:gauge} Assume $g$ is of class $\smash{C_\loc^{k,\alpha}}$, where $k\in\N_0$ and $\alpha\in (0,1)$. Given $\smash{x_0\in M}$, let $r$ and $\Phi$ be the $p$-harmonic coordinate radius and map from \cref{thm:coord} and set $\smash{\widehat g := \Phi_*g}$ on the open set $\Phi(B_r(x_0))$. Then  $\smash{\widehat g^{\mu\mu}<0}$ on $\Phi(B_r(x_0))$ and for every $\smash{\varphi\in C_\comp^\infty(\Phi(B_r(x_0)))}$,
\begin{align*}
    \int_{\Phi(B_r(x_0))}\big\vert\!\det\widehat g\big\vert^{1/2}\,\big[{-\widehat g^{\mu\mu}}\big]^{(p-2)/2}\,\widehat g^{\alpha\mu}\,\partial_\alpha\varphi\d\Leb^{1+d} = 0,
\end{align*}
without summation over $\mu$.
\end{corollary}

\begin{proof} We time-orient $\smash{\widehat{g}}$ by the push-forward of the time-orientation $X$ under $\Phi$. Then by \cref{thm:coord}, $\smash{\Phi^{-1}}$ is a $C^1$-diffeomorphism from $\Phi(B_r(x_0))$ onto $B_r(x_0)$; by construction, it is a time-orientation-preserving isometry between $\smash{\widehat g}$ and $g$ at every point of $\Phi(B_r(x_0))$. By \cref{Le:Invarp}, the coordinate function $\smash{y^\mu := \tau^\mu\circ\Phi^{-1}}$ is a $p$-harmonic temporal function for $\smash{\widehat g}$ on $\Phi(B_r(x_0))$. Its differential has the constant coefficients $\smash{\delta^\mu_\beta}$, so that $\smash{\widehat g^{\alpha\beta}\,\partial_\alpha y^\mu\,\partial_\beta y^\mu = \widehat g^{\mu\mu}}$ and $\smash{\widehat g^{\alpha\beta}\,\partial_\beta y^\mu = \widehat g^{\alpha\mu}}$. As $\smash{\widehat g^{\mu\mu}}$ is continuous, the claim follows from  \eqref{Eq:Weak}.
\end{proof}

\subsection{Proofs of \cref{thm:iso,thm:bilip}} Most work is invested in the proof of \cref{thm:bilip}, which includes \cref{thm:iso} as a special case.

\begin{proof}[Proof of \cref{thm:bilip}, $\alpha\in (0,1)$] As the claim is local, we may and will assume $\psi$ is a homeomorphism. Let $\alpha\in(0,1)$ and $x_0\in M$. By \cref{thm:coord} applied to $g'$, there are $r>0$ and a $\smash{C^{k+1,\alpha}}$-diffeomorphism $\smash{\Phi'}$ from $\smash{B_r(\psi(x_0))}$ onto its image whose components are $p$-harmonic temporal functions for $g'$. Note $\smash{\psi^{-1}(B_r(\psi(x_0)))}$ is an open neighborhood of $x_0$. Therefore, the restriction of $\psi$ to $\smash{\psi^{-1}(B_r(\psi(x_0)))}$ is a locally bi-Lipschitz homeomorphism onto $B_r(\psi(x_0))$ which is a time-orientation-preserving isometry on $\psi^{-1}(B_r(\psi(x_0)))$. By \cref{Le:Invarp}, the components of $\Phi'\circ\psi$ are $p$-harmonic temporal functions for $g$ on $\smash{\psi^{-1}(B_r(\psi(x_0)))}$. By \cref{thm:reg}, they are of class $\smash{C^{k+1,\alpha}_\loc}$ on $\smash{\psi^{-1}(B_r(\psi(x_0)))}$. This shows $\smash{\psi = \Phi'^{-1}\circ(\Phi'\circ\psi)}$ is of class $\smash{C^{k+1,\alpha}_\loc}$ on $\smash{\psi^{-1}(B_r(\psi(x_0)))}$. As $x_0$ was arbitrary, this proves the claim.
\end{proof}

To treat the endpoint of the Hölder scale, we need the following calculus rules. The proof is straightforward and follows from the standard chain and Leibniz rules for Sobolev functions, cf.~e.g.~Evans--Gariepy \cite{evans-gariepy1992}*{§4}.

\begin{lemma}[Calculus rules]\label{lem:H2calc} Let $\smash{V,W\subset\R^{1+d}}$ be open sets. Let $\Xi\colon V\to W$ be a $C^1$-diffeomorphism whose components belong to $\smash{H^2_\loc(V)}$. Then the following hold.
\begin{enumerate}[label=\textnormal{(\roman*)}]
    \item\label{La:C1} The components of $\smash{\Xi^{-1}}$ belong to $\smash{H^2_\loc(W)}$.
    \item\label{La:C2} If $\smash{f\in C^1(W)\cap H^2_\loc(W)}$, we have $\smash{f\circ\Xi\in H^2_\loc(V)}$.
\end{enumerate}
\end{lemma}

\begin{proof}[Proof of \cref{thm:bilip}, $\alpha=1$] As in the first part of the proof, we may and will assume $\psi$ is a homeomorphism. We fix $\sigma\in(0,1)$ and proceed in four steps.

\emph{Step 1. Preliminary regularity.} Since $g$ and $g'$ are of class $\smash{C^{k,\sigma}_\loc}$, the part of the proof above shows $\psi$ is of class $\smash{C^{k+1,\sigma}_\loc}$. In particular, $\psi$ is of class $C^1$. Hence, both sides of the $\smash{\Leb^{1+d}}$-a.e.~isometry of $\psi$ depend continuously on $x\in M$. As they coincide for $\smash{\Leb^{1+d}}$-a.e.~$x\in M$, they coincide for every $x\in M$. In coordinates, 
\begin{align}\label{Eq:PullbackCoord}
    g_{\mu\nu} = g'_{\gamma\delta}\circ\psi\,\partial_\mu\psi^\gamma\,\partial_\nu\psi^\delta
\end{align}
everywhere on $M$. Since $g$ and $g'$ are nondegenerate, taking determinants implies $\rmd\psi(x)$ is invertible for every $x\in M$. As $\psi$ is a homeomorphism, the inverse function theorem shows $\psi$ is a $C^1$-diffeomorphism from $M$ onto $M'$.

\emph{Step 2. Second derivatives for locally Lipschitz metrics.} Let $k=0$. We claim the components of $\psi$ belong to $\smash{H^2_\loc(M)}$. Let $x_0\in M$. Moreover, let $r$ and $\Phi'$ be as in the first part of the proof, with $\sigma$ in place of $\alpha$. Since the metric $g'$ is locally Lipschitz, \cref{lem:H2} shows the components of $\Phi'$ belong to $\smash{H^2_\loc(B_r(\psi(x_0)))}$. As $\Phi'$ is a $C^1$-diffeomorphism onto its image, \cref{lem:H2calc} implies the components of $\smash{\Phi'^{-1}}$ belong to $\smash{H^2_\loc(\Phi'(B_r(\psi(x_0))))}$; moreover, they are of class $C^1$, as $\Phi'$ is a $C^1$-diffeomorphism. Next, let us consider the map $\Xi:=\Phi'\circ\psi$ on the open set $\psi^{-1}(B_r(\psi(x_0)))$. By Step 1, it is a $C^1$-diffeomorphism from $\psi^{-1}(B_r(\psi(x_0)))$ onto $\Phi'(B_r(\psi(x_0)))$. Then by \cref{Le:Invarp}, its components are $p$-harmonic temporal functions for $g$ on $\psi^{-1}(B_r(\psi(x_0)))$. Since $g$ is locally Lipschitz, \cref{lem:H2} shows the components of $\Xi$ belong to $\smash{H^2_\loc(\psi^{-1}(B_r(\psi(x_0))))}$. Applying \cref{lem:H2calc} to all components of $\smash{\Phi'^{-1}}$ in place of $f$, the components of the composition $\smash{\psi = \Phi'^{-1}\circ\Xi}$ belong to $\smash{H^2_\loc(\psi^{-1}(B_r(\psi(x_0))))}$. Since $x_0$ was arbitrary, the claim follows.

\emph{Step 3. Lipschitz continuity of $\rmd\psi$ for locally Lipschitz metrics.} Retain the assumption $k=0$. We now differentiate \eqref{Eq:PullbackCoord}. The left-hand side is locally Lipschitz. On the other hand, the function $\smash{g'_{\gamma\delta}\circ\psi}$ in the right-hand side is locally Lipschitz as a composition of a locally Lipschitz function and a map of class $C^1$, cf.~Step 1. Since $\smash{\psi^{-1}}$ is locally Lipschitz, it maps negligible sets to negligible sets with respect to $\smash{\Leb^{1+d}}$. Hence, by Rademacher's theorem, $\smash{g'_{\gamma\delta}}$ is differentiable at $\psi(x)$ for $\smash{\Leb^{1+d}}$-a.e.~$x\in M$; at every such point the chain rule gives
\begin{align*}
    \partial_\kappa\big[g'_{\gamma\delta}\circ\psi\big] = \partial_\theta g'_{\gamma\delta}\circ\psi\,\partial_\kappa\psi^\theta.
\end{align*}
The first derivatives of $\psi$ are continuous and belong to $\smash{H^1_\loc}$ by Step 2. Thus, the product rule applies to the right-hand side of \eqref{Eq:PullbackCoord}. We abbreviate
\begin{align*}
    T_{\kappa\mu\nu} &:= \partial_\theta g'_{\gamma\delta}\circ\psi\,\partial_\kappa\psi^\theta\,\partial_\mu\psi^\gamma\,\partial_\nu\psi^\delta,\\
    S_{\kappa\mu\nu} &:= g'_{\gamma\delta}\circ\psi\,\partial_\kappa\partial_\mu\psi^\gamma\,\partial_\nu\psi^\delta.
\end{align*}
Using the symmetry of $g'$ in the last term, differentiating \eqref{Eq:PullbackCoord} yields
\begin{align*}
    \partial_\kappa g_{\mu\nu} &= T_{\kappa\mu\nu} + g'_{\gamma\delta}\circ\psi\,\partial_\kappa\partial_\mu\psi^\gamma\,\partial_\nu\psi^\delta + g'_{\gamma\delta}\circ\psi\,\partial_\mu\psi^\gamma\,\partial_\kappa\partial_\nu\psi^\delta\\
    &= T_{\kappa\mu\nu} + S_{\kappa\mu\nu} + S_{\kappa\nu\mu}\quad\Leb^{1+d}\textnormal{-a.e.}
\end{align*}
We add this identity to the one obtained by exchanging $\kappa$ and $\mu$ and subtract the one obtained by exchanging $\kappa$ and $\nu$. This gives
\begin{align*}
    &\partial_\kappa g_{\mu\nu} + \partial_\mu g_{\kappa\nu} - \partial_\nu g_{\mu\kappa} - T_{\kappa\mu\nu} - T_{\mu\kappa\nu} + T_{\nu\mu\kappa}\\
    &\qquad\qquad = S_{\kappa\mu\nu} + S_{\kappa\nu\mu} + S_{\mu\kappa\nu} + S_{\mu\nu\kappa} - S_{\nu\mu\kappa} - S_{\nu\kappa\mu}\\
    &\qquad\qquad = 2\,S_{\kappa\mu\nu}\quad\Leb^{1+d}\textnormal{-a.e.};
\end{align*}
in the last step, we employed that $S_{\kappa\mu\nu}$ is symmetric in its first two indices. The first derivatives of $g$ belong to $\smash{L^\infty_\loc}$ since $g$ is locally Lipschitz. Moreover, the functions $\smash{\partial_\theta g'_{\gamma\delta}\circ\psi}$ belong to $\smash{L^\infty_\loc}$ as $g'$ is locally Lipschitz and $\psi$ is a homeomorphism. In addition, the first derivatives of $\psi$ are continuous. Consequently, the left-hand side of the preceding identity belongs to $\smash{L^\infty_\loc}$ and so does the expression $S_{\kappa\mu\nu}$.

It remains to solve for the second derivatives of $\psi$. We abbreviate
\begin{align*}
    N_{\gamma\nu} := g'_{\gamma\delta}\circ\psi\,\partial_\nu\psi^\delta,
\end{align*}
which is continuous. Since $g'$ is nondegenerate and $\rmd\psi$ is invertible by Step 1, its values are invertible matrices. Let $\smash{N^{\nu\theta}}$ denote the entries of its (not necessarily symmetric) inverse, which are continuous. By definition, $\smash{S_{\kappa\mu\nu} = N_{\gamma\nu}\,\partial_\kappa\partial_\mu\psi^\gamma}$, so
\begin{align*}
    \partial_\kappa\partial_\mu\psi^\theta = S_{\kappa\mu\nu}\,N^{\nu\theta}\quad\Leb^{1+d}\textnormal{-a.e.}
\end{align*}
The right-hand side belongs to $\smash{L^\infty_\loc(M)}$. Thus, the first derivatives of $\psi$ belong to $\smash{W^{1,\infty}_\loc(M)}$, whence they are locally Lipschitz. This shows $\psi$ is of class $\smash{C^{1,1}_\loc(M)}$.

\emph{Step 4. The case $k\in\N$.} Now assume $k\in\N$. By Step 1, the map $\psi$ is of class $\smash{C^{k+1,\sigma}_\loc}$. In particular, it is of regularity $C^2$ and $g$ and $g'$ are of class $C^1$. Thus, the computations from Step 3 hold in the classical sense at every point of $M$. With the notation from this step, this entails the following pointwise identity on $M$:
\begin{align*}
    \partial_\kappa\partial_\mu\psi^\theta = \frac{1}{2}\,\big[\partial_\kappa g_{\mu\nu} + \partial_\mu g_{\kappa\nu} - \partial_\nu g_{\mu\kappa} - T_{\kappa\mu\nu} - T_{\mu\kappa\nu} + T_{\nu\mu\kappa}\big]\,N^{\nu\theta}.
\end{align*}
This is the transformation rule of the Christoffel symbols, cf.~Calabi--Hartman \cite{calabi-hartman1970}, (5.2) therein. The first derivatives of $g$ and $g'$ are of class $\smash{C^{k-1,1}_\loc}$; $\psi$ and its first derivatives are of class $\smash{C^{k,\sigma}_\loc}$, hence of class $\smash{C^{k-1,1}_\loc}$. Since compositions, products, and inverses of maps of class $\smash{C^{k-1,1}_\loc}$ are again of this class, the right-hand side of the previous identity is of class $\smash{C^{k-1,1}_\loc}$. Therefore, $\psi$ is of class $\smash{C^{k+1,1}_\loc}$.
\end{proof}

\begin{proof}[Proof of \cref{thm:iso}] In view of \cref{thm:bilip} established above, we only have to prove $\psi$ is time-orientation-preserving up to a change of sign in the underlying time-orientation of $g'$, which does not affect the hypotheses on $g'$.

Let $x_0\in M$. As $g$ and $g'$ are nondegenerate, taking determinants in the identity $\smash{\psi^*g'(x_0)=g(x_0)}$ shows that $\rmd\psi(x_0)$ is invertible. By the inverse function theorem, there is a connected open neighborhood $U\subset M$ of $x_0$ such that $\psi$ restricts to a $C^1$-diffeomorphism from $U$ onto the open set $\psi(U)$; this restriction is a locally bi-Lipschitz homeomorphism. For every $x\in U$, $\rmd\psi(x)X(x)$ and $X'\circ\psi(x)$ are timelike with respect to $g'\circ\psi(x)$. By the reverse Cauchy--Schwarz inequality, this implies the quantity $g'\circ\psi(x)(\rmd \psi(x)X(x),X'\circ\psi(x))$, which depends continuously on $x\in U$, is nowhere vanishing on $U$. Since $U$ is connected, it has a sign. If it is positive, we replace $X'$ by $-X'$. In either case, $\psi$ is a time-orientation-preserving isometry at every point of $U$. \cref{thm:bilip} readily shows that $\psi$ is of class $\smash{C^{k+1,\alpha}_\loc}$ on $U$. The arbitrariness of $x_0$ implies the claim.
\end{proof}

\begin{corollary}[Uniqueness of Hölder structures]\label{cor:atlas} Let $k\in\N_0$ and $\alpha\in(0,1)$. Let $M$ be a $C^1$-manifold of dimension $1+d$ with a continuous Lorentzian metric $g$. Then there is at most one $\smash{C^{k+1,\alpha}}$-structure on $M$ which is compatible with the given $C^1$-structure and in whose charts $g$ has coefficients of class $\smash{C^{k,\alpha}_\loc}$. 

If it exists, the $p$-harmonic coordinates from \cref{thm:coord} form an atlas of it.
\end{corollary}

\begin{proof} Given two such structures, every transition map between a chart of the first and a chart of the second structure is of class $C^1$ and pulls the coefficients of $g$ in the second chart back to those in the first chart. By \cref{thm:iso}, it is of class $\smash{C^{k+1,\alpha}_\loc}$ and so is its inverse. Therefore, the two structures coincide. 

The last claim follows by \cref{thm:coord} applied in the charts of this structure.
\end{proof}

\subsection{Proof of \cref{thm:tau}} Let $g$ be a continuous Lorentzian metric on $M$. A curve $\gamma\colon[0,1]\to M$ is called \textbf{future-directed causal} if it is locally Lipschitz and $\smash{\dot\gamma_s\in T_{\gamma_s}M}$ is a future-directed causal vector for $\smash{\Leb^1}$-a.e.~$s\in[0,1]$. Its \textbf{length} is
\begin{align*}
    \Len_g(\gamma) := \int_{[0,1]}\vert\dot\gamma_s\vert_{g(\gamma_s)}\d s.
\end{align*}
The \textbf{time separation function} $l\colon M\times M\to\R\cup\{-\infty,\infty\}$ of $g$ is defined by letting $l(x,y)$ be the supremum of $\Len_g(\gamma)$ over all causal curves $\gamma$ from $x$ to $y$, with the usual convention $\sup\emptyset:=-\infty$. We call $g$ \textbf{strongly causal} if for every point $x\in M$ and every open neighborhood $U\subset M$ of $x$, there exists an open neighborhood $V\subset U$ of $x$ such that every causal curve with endpoints in $V$ is contained in $U$; we refer to Minguzzi \cite{minguzzi2019-causality} for an overview on causality theory. 

\begin{definition}[Time-separation-preserving map]\label{Def:TimeSepPres} We say a map $\psi\colon M\to M'$ preserves the time separation functions $l$ and $l'$ of two continuous Lorentzian metrics $g$ and $g'$ on $M$ and $M'$, respectively, if every $x,y\in M$ satisfy 
\begin{align*}
l_+'(\psi(x),\psi(y)) = l_+(x,y).    
\end{align*}
\end{definition}

The time separation function determines the metric to first order.  The proof of Ryborz \cite[Lem.~4.13]{ryborz2026+} carries over  verbatim, with the vector path below in place of the difference quotients of the causal path considered therein. 

\begin{lemma}[First-order expansion]\label{lem:tau} Let $g$ be a continuous Lorentzian metric on $M$ which is strongly causal and let $x\in M$. Let $\varepsilon>0$ and let $v\colon(-\varepsilon,\varepsilon)\to \R^{1+d}$ be continuous at $0$. If $v_0$ is a future-directed timelike vector at $x$,
\begin{align*}
    \lim_{t\to 0^+}\frac{l_+(x,x+tv_t)}{t} = \big\vert v_0\big\vert_{g(x)};
\end{align*}
otherwise, this limit is $0$.
\end{lemma}

\begin{proof}[Proof of \cref{thm:tau}] By \cref{thm:bilip}, it suffices to show $\psi$ is a time-orientation preserving isometry at every $x\in M$ where it is differentiable. Since $\psi$ is locally bi-Lipschitz, $\rmd\psi(x)$ is invertible. Let $w\in T_xM$ be a future-directed $g(x)$-timelike vector. For sufficiently small $\varepsilon>0$, define $v\colon (-\varepsilon,\varepsilon)\to \R^{1+d}$ by
\begin{align*} 
    v_t := \frac{\psi(x+tw)-\psi(x)}{t},
\end{align*}
so that $\psi(x+tw)=\psi(x) +tv_t$ for every $t\in (-\varepsilon,\varepsilon)$ and $v_t\to \rmd\psi(x)w$ as $t\to 0^+$. Since $\psi$ preserves the time separation functions, \cref{lem:tau} for $g$ gives
\begin{align*}
    \lim_{t\to 0^+}\frac{l_+'(\psi(x),\psi(x)+t v_t)}{t} &= \lim_{t\to 0^+} \frac{l_+'(\psi(x),\psi(x + tw))}{t}\\
    &=\lim_{t\to 0^+}\frac{l_+(x,x+tw)}{t}\\
    &= \big\vert w\big\vert_{g(x)}.
\end{align*}
Since the right-hand side is positive, by \cref{lem:tau} for $g'$ this forces $\rmd\psi(x)w$ to be a future-directed $g'\circ\psi(x)$-timelike vector at $\psi(x)$ with 
\begin{align*}
    \vert \rmd\psi(x)w\vert_{g'\circ\psi(x)} = \vert w\vert_{g(x)}
\end{align*}
Hence, $\psi$ is a time-orientation-preserving map at $x$ and the two quadratic forms $\psi^*g'(x)$ and $g(x)$ coincide on the nonempty open set of all future-directed $g(x)$-timelike vectors at $x$. As both quadratic forms are polynomials and future-directed $g(x)$-timelike vectors form a basis of $T_xM$, we obtain $\psi^*g'(x)=g(x)$ as quadratic forms on $\smash{\R^{1+d}}$, whence $\psi$ is an isometry at $x$. Then the arbitrariness of $x$ readily establishes the claim.
\end{proof}

\addtocontents{toc}{\protect\setcounter{tocdepth}{1}}

\subsection*{Declaration on the use of AI} During the preparation of this work, the author made use of the large language model Claude (Anthropic). The model produced a first draft of the manuscript, including drafts of the proofs, on the basis of the author's ideas and instructions to construct $p$-harmonic coordinates and apply them to the regularity of spacetime isometries. It was further used to check arguments and to locate and verify references. The author subsequently rewrote the entire text, verified every statement and proof line by line, and takes full responsibility for the content of this article.

\subsection*{Conflict of interest} The author declares no conflict of interest.

\subsection*{Data availability} Data sharing is not applicable to this article, as no datasets were generated or analyzed.

\bibliographystyle{amsrefs}
\bibliography{library}

@Article{mccann2020,
  author   = {McCann, Robert J.},
  journal  = {Camb. J. Math.},
  title    = {Displacement convexity of {B}oltzmann's entropy characterizes the strong energy condition from general relativity},
  year     = {2020},
  issn     = {2168-0930,2168-0949},
  number   = {3},
  pages    = {609--681},
  volume   = {8},
  doi      = {10.4310/CJM.2020.v8.n3.a4},
  fjournal = {Cambridge Journal of Mathematics},
  mrclass  = {53C50 (49Q22 53C21 58Z05 82C35 83C99)},
  mrnumber = {4192570},
  url      = {https://doi.org/10.4310/CJM.2020.v8.n3.a4},
}

@Article{mondino-suhr2022,
  author   = {Mondino, Andrea and Suhr, Stefan},
  journal  = {J. Eur. Math. Soc. (JEMS)},
  title    = {An optimal transport formulation of the {E}instein equations of general relativity},
  year     = {2023},
  issn     = {1435-9855,1435-9863},
  number   = {3},
  pages    = {933--994},
  volume   = {25},
  doi      = {10.4171/jems/1188},
  fjournal = {Journal of the European Mathematical Society (JEMS)},
  mrclass  = {53C80 (49Q22 53C50 82C35 83C05)},
  mrnumber = {4577957},
  url      = {https://doi.org/10.4171/jems/1188},
}

@Misc{beran-braun-calisti-gigli-mccann-ohanyan-rott-samann+-,
  author = {Beran, Tobias and Braun, Mathias and Calisti, Matteo and Gigli, Nicola and McCann, Robert J. and Ohanyan, Argam and Rott, Felix and S\"amann, Clemens},
  title  = {{A} nonlinear d'{A}lembert comparison theorem and causal differential calculus on metric measure spacetimes},
  year   = {Preprint, arXiv:2408.15968},
}

@Misc{dephilippis-gigli-mccann-ohanyan-picerni-zanardini+,
  author = {De Philippis, Guido and Gigli, Nicola and McCann, Robert J. and Ohanyan, Argam and Picerni, Marco and Zanardini, Matteo},
  year   = {In preparation},
}

@Article{minguzzi2019-causality,
  author    = {Minguzzi, Ettore},
  journal   = {Living Reviews in Relativity},
  title     = {Lorentzian causality theory},
  year      = {2019},
  issn      = {1433-8351},
  number    = {3},
  pages     = {202 pp.},
  volume    = {22},
  doi       = {10.1007/s41114-019-0019-x},
  publisher = {Springer Science and Business Media LLC},
  url       = {http://dx.doi.org/10.1007/s41114-019-0019-x},
}

@article{deturck-kazdan1981,
    AUTHOR = {DeTurck, Dennis M. and Kazdan, Jerry L.},
     TITLE = {Some regularity theorems in {R}iemannian geometry},
   JOURNAL = {Ann. Sci. \'Ecole Norm. Sup. (4)},
  FJOURNAL = {Annales Scientifiques de l'\'Ecole Normale Sup\'erieure.
              Quatri\`eme S\'erie},
    VOLUME = {14},
      YEAR = {1981},
    NUMBER = {3},
     PAGES = {249--260},
      ISSN = {0012-9593},
   MRCLASS = {53B21},
  MRNUMBER = {644518},
MRREVIEWER = {V.\ W\"unsch},
       URL = {http://www.numdam.org/item?id=ASENS_1981_4_14_3_249_0},
}

@article{caponio-ohanyan-ohta2024+,
    AUTHOR = {Caponio, Erasmo and Ohanyan, Argam and Ohta, Shinichi},
     TITLE = {Splitting theorems for weighted {F}insler spacetimes via the
              {$p$}-d'{A}lembertian: beyond the {B}erwald case},
   JOURNAL = {J. Reine Angew. Math.},
  FJOURNAL = {Journal f\"ur die Reine und Angewandte Mathematik. [Crelle's
              Journal]},
    VOLUME = {838},
      YEAR = {2026},
     PAGES = {221--253},
      ISSN = {0075-4102,1435-5345},
   MRCLASS = {53C60 (53C24 53C50)},
  MRNUMBER = {5124089},
       DOI = {10.1515/crelle-2026-0045},
       URL = {https://doi.org/10.1515/crelle-2026-0045},
}

@article{hartman1952,
    AUTHOR = {Hartman, Philip},
     TITLE = {On unsmooth two-dimensional {R}iemannian metrics},
   JOURNAL = {Amer. J. Math.},
  FJOURNAL = {American Journal of Mathematics},
    VOLUME = {74},
      YEAR = {1952},
     PAGES = {215--226},
      ISSN = {0002-9327,1080-6377},
   MRCLASS = {53.0X},
  MRNUMBER = {47379},
MRREVIEWER = {A.\ Fialkow},
       DOI = {10.2307/2372079},
       URL = {https://doi.org/10.2307/2372079},
}

@article{sabitov1993,
  author    = {Sabitov, I. Kh.},
  title     = {To the question of smoothness of isometries},
  journal   = {Siberian Math. J.},
  volume    = {34},
  number    = {4},
  year      = {1993},
  pages     = {741--748},
  doi       = {10.1007/BF00975178},
  note      = {Translated from Sibirsk. Mat. Zh. \textbf{34} (1993), no. 4, 169--176}
}

@article{shefel1982,
    AUTHOR = {Shefel{'}, S. Z.},
     TITLE = {Smoothness of a conformal mapping of {R}iemannian spaces},
   JOURNAL = {Sibirsk. Mat. Zh.},
  FJOURNAL = {Akademiya Nauk SSSR. Sibirskoe Otdelenie. Sibirski\u i\
              Matematicheski\u i\ Zhurnal},
    VOLUME = {23},
      YEAR = {1982},
    NUMBER = {1},
     PAGES = {153--159, 222},
      ISSN = {0037-4474},
   MRCLASS = {53C20},
  MRNUMBER = {651886},
MRREVIEWER = {A.\ V.\ Pogorelov},
}

@article{liimatainen-salo2014,
    AUTHOR = {Liimatainen, Tony and Salo, Mikko},
     TITLE = {{$n$}-harmonic coordinates and the regularity of conformal
              mappings},
   JOURNAL = {Math. Res. Lett.},
  FJOURNAL = {Mathematical Research Letters},
    VOLUME = {21},
      YEAR = {2014},
    NUMBER = {2},
     PAGES = {341--361},
      ISSN = {1073-2780,1945-001X},
   MRCLASS = {53C20 (35B65 35J92 53A30)},
  MRNUMBER = {3247061},
MRREVIEWER = {Zhi\ Ren\ Jin},
       DOI = {10.4310/MRL.2014.v21.n2.a11},
       URL = {https://doi.org/10.4310/MRL.2014.v21.n2.a11},
}

@article{colombini-degiorgi-spagnolo1979,
    AUTHOR = {Colombini, Ferruccio and De Giorgi, Ennio and Spagnolo,
              Sergio},
     TITLE = {Sur les \'equations hyperboliques avec des coefficients qui ne
              d\'ependent que du temps},
   JOURNAL = {Ann. Scuola Norm. Sup. Pisa Cl. Sci. (4)},
  FJOURNAL = {Annali della Scuola Normale Superiore di Pisa. Classe di
              Scienze. Serie IV},
    VOLUME = {6},
      YEAR = {1979},
    NUMBER = {3},
     PAGES = {511--559},
      ISSN = {0391-173X,2036-2145},
   MRCLASS = {35L15},
  MRNUMBER = {553796},
MRREVIEWER = {P.\ Ungar},
       URL = {http://www.numdam.org/item?id=ASNSP_1979_4_6_3_511_0},
}

@article{klainerman-rodnianski2005,
    AUTHOR = {Klainerman, Sergiu and Rodnianski, Igor},
     TITLE = {Rough solutions of the {E}instein-vacuum equations},
   JOURNAL = {Ann. of Math. (2)},
  FJOURNAL = {Annals of Mathematics. Second Series},
    VOLUME = {161},
      YEAR = {2005},
    NUMBER = {3},
     PAGES = {1143--1193},
      ISSN = {0003-486X,1939-8980},
   MRCLASS = {58J45 (35Q75)},
  MRNUMBER = {2180400},
MRREVIEWER = {Lars\ \AA ke\ Andersson},
       DOI = {10.4007/annals.2005.161.1143},
       URL = {https://doi.org/10.4007/annals.2005.161.1143},
}

@article{lanczos1922,
  author  = {Lanczos, Cornelius},
  title   = {Ein vereinfachendes {K}oordinatensystem f{\"u}r die {E}insteinschen {G}ravitationsgleichungen},
  journal = {Phys. Z.},
  volume  = {23},
  year    = {1922},
  pages   = {537--539}
}

@article{liimatainen-salo2016,
    AUTHOR = {Liimatainen, Tony and Salo, Mikko},
     TITLE = {Local gauge conditions for ellipticity in conformal geometry},
   JOURNAL = {Int. Math. Res. Not. IMRN},
  FJOURNAL = {International Mathematics Research Notices. IMRN},
      YEAR = {2016},
    NUMBER = {13},
     PAGES = {4058--4077},
      ISSN = {1073-7928,1687-0247},
   MRCLASS = {53A30 (58J05)},
  MRNUMBER = {3544628},
MRREVIEWER = {Gang\ Li},
       DOI = {10.1093/imrn/rnv255},
       URL = {https://doi.org/10.1093/imrn/rnv255},
}

@book{dedonder1921,
  author    = {De Donder, Th{\'e}ophile},
  title     = {La gravifique einsteinienne},
  publisher = {Gauthier-Villars},
  address   = {Paris},
  year      = {1921}
}

@article{einstein1916-integration,
  author  = {Einstein, Albert},
  title   = {N{\"a}herungsweise {I}ntegration der {F}eldgleichungen der {G}ravitation},
  journal = {Sitzungsber. K{\"o}nigl. Preuss. Akad. Wiss. Berlin},
  year    = {1916},
  pages   = {688--696},
  url     = {https://einsteinpapers.press.princeton.edu/vol6-doc/375}
}

@article{julin-liimatainen-salo2017,
    AUTHOR = {Julin, Vesa and Liimatainen, Tony and Salo, Mikko},
     TITLE = {{$p$}-harmonic coordinates for {H}\"older metrics and
              applications},
   JOURNAL = {Comm. Anal. Geom.},
  FJOURNAL = {Communications in Analysis and Geometry},
    VOLUME = {25},
      YEAR = {2017},
    NUMBER = {2},
     PAGES = {395--430},
      ISSN = {1019-8385,1944-9992},
   MRCLASS = {53A30 (31C12 47J05 58C25)},
  MRNUMBER = {3690246},
MRREVIEWER = {Jui-Tang\ Chen},
       DOI = {10.4310/CAG.2017.v25.n2.a5},
       URL = {https://doi.org/10.4310/CAG.2017.v25.n2.a5},
}

@article{hartman-wintner1953,
    AUTHOR = {Hartman, Philip and Wintner, Aurel},
     TITLE = {On the existence of {R}iemannian manifolds which cannot carry
              non-constant analytic or harmonic functions in the small},
   JOURNAL = {Amer. J. Math.},
  FJOURNAL = {American Journal of Mathematics},
    VOLUME = {75},
      YEAR = {1953},
     PAGES = {260--276},
      ISSN = {0002-9327,1080-6377},
   MRCLASS = {53.0X},
  MRNUMBER = {55016},
MRREVIEWER = {A.\ Lichnerowicz},
       DOI = {10.2307/2372453},
       URL = {https://doi.org/10.2307/2372453},
}

@article{matveev-troyanov2017,
    AUTHOR = {Matveev, Vladimir S. and Troyanov, Marc},
     TITLE = {The {M}yers-{S}teenrod theorem for {F}insler manifolds of low
              regularity},
   JOURNAL = {Proc. Amer. Math. Soc.},
  FJOURNAL = {Proceedings of the American Mathematical Society},
    VOLUME = {145},
      YEAR = {2017},
    NUMBER = {6},
     PAGES = {2699--2712},
      ISSN = {0002-9939,1088-6826},
   MRCLASS = {53B40 (35B65 53C60)},
  MRNUMBER = {3626522},
MRREVIEWER = {Behroz\ Bidabad},
       DOI = {10.1090/proc/13407},
       URL = {https://doi.org/10.1090/proc/13407},
}

@article{wintner1952,
    AUTHOR = {Wintner, Aurel},
     TITLE = {On isometric surfaces},
   JOURNAL = {Amer. J. Math.},
  FJOURNAL = {American Journal of Mathematics},
    VOLUME = {74},
      YEAR = {1952},
     PAGES = {198--214},
      ISSN = {0002-9327,1080-6377},
   MRCLASS = {52.0X},
  MRNUMBER = {47351},
MRREVIEWER = {H.\ Busemann},
       DOI = {10.2307/2372078},
       URL = {https://doi.org/10.2307/2372078},
}

@article{taylor2006,
    AUTHOR = {Taylor, Michael},
     TITLE = {Existence and regularity of isometries},
   JOURNAL = {Trans. Amer. Math. Soc.},
  FJOURNAL = {Transactions of the American Mathematical Society},
    VOLUME = {358},
      YEAR = {2006},
    NUMBER = {6},
     PAGES = {2415--2423},
      ISSN = {0002-9947,1088-6850},
   MRCLASS = {53C21 (35B65 53C42)},
  MRNUMBER = {2204038},
MRREVIEWER = {Federico\ S\'anchez-Bringas},
       DOI = {10.1090/S0002-9947-06-04090-6},
       URL = {https://doi.org/10.1090/S0002-9947-06-04090-6},
}

@article{reshetnyak1978,
    AUTHOR = {Reshetnjak, Ju. G.},
     TITLE = {Differential properties of quasiconformal mappings and
              conformal mappings of {R}iemannian spaces},
   JOURNAL = {Sibirsk. Mat. Zh.},
  FJOURNAL = {Akademiya Nauk SSSR. Sibirskoe Otdelenie. Sibirski\u i\
              Matematicheski\u i\ Zhurnal},
    VOLUME = {19},
      YEAR = {1978},
    NUMBER = {5},
     PAGES = {1166--1183, 1216},
      ISSN = {0037-4474},
   MRCLASS = {53B20 (30C60)},
  MRNUMBER = {508507},
MRREVIEWER = {Valentin\ Boju},
}

@article{lytchak-yaman2006,
    AUTHOR = {Lytchak, Alexander and Yaman, Asli},
     TITLE = {On {H}\"older continuous {R}iemannian and {F}insler metrics},
   JOURNAL = {Trans. Amer. Math. Soc.},
  FJOURNAL = {Transactions of the American Mathematical Society},
    VOLUME = {358},
      YEAR = {2006},
    NUMBER = {7},
     PAGES = {2917--2926},
      ISSN = {0002-9947,1088-6850},
   MRCLASS = {53C22 (53C20 53C60)},
  MRNUMBER = {2216252},
MRREVIEWER = {Igor\ V.\ Nikolaev},
       DOI = {10.1090/S0002-9947-06-04195-X},
       URL = {https://doi.org/10.1090/S0002-9947-06-04195-X},
}

@article{hartman1958,
    AUTHOR = {Hartman, Philip},
     TITLE = {On isometries and on a theorem of {L}iouville},
   JOURNAL = {Math. Z.},
  FJOURNAL = {Mathematische Zeitschrift},
    VOLUME = {69},
      YEAR = {1958},
     PAGES = {202--210},
      ISSN = {0025-5874,1432-1823},
   MRCLASS = {53.00 (35.00)},
  MRNUMBER = {108809},
MRREVIEWER = {A.\ Douglis},
       DOI = {10.1007/BF01187401},
       URL = {https://doi.org/10.1007/BF01187401},
}

@Misc{braun+,
  author={Braun, Mathias},
  title={},
  year ={In preparation},
}

@Article{eschenburg1988,
  author     = {Eschenburg, J.-H.},
  journal    = {J. Differential Geom.},
  title      = {The splitting theorem for space-times with strong energy condition},
  year       = {1988},
  issn       = {0022-040X,1945-743X},
  number     = {3},
  pages      = {477--491},
  volume     = {27},
  fjournal   = {Journal of Differential Geometry},
  mrclass    = {53C50},
  mrnumber   = {940115},
  mrreviewer = {Paul\ E.\ Ehrlich},
  url        = {http://projecteuclid.org/euclid.jdg/1214442005},
}

@Article{bartnik1988,
  author     = {Bartnik, Robert},
  journal    = {Comm. Math. Phys.},
  title      = {Remarks on cosmological spacetimes and constant mean curvature surfaces},
  year       = {1988},
  issn       = {0010-3616,1432-0916},
  number     = {4},
  pages      = {615--624},
  volume     = {117},
  fjournal   = {Communications in Mathematical Physics},
  mrclass    = {53C50 (53C42 83C99)},
  mrnumber   = {953823},
  mrreviewer = {Judith\ M.\ Arms},
  url        = {http://projecteuclid.org/euclid.cmp/1104161820},
}

@article{galloway-mccann-ohanyan2026+,
    AUTHOR = {Galloway, Gregory J. and McCann, Robert J. and Ohanyan, Argam},
     TITLE = {Failure of local equi-{L}ipschitzness for families of
              {L}orentz distances to {C}auchy surface foliations},
   JOURNAL = {Gen. Relativity Gravitation},
  FJOURNAL = {General Relativity and Gravitation},
    VOLUME = {58},
      YEAR = {2026},
    NUMBER = {9},
     PAGES = {Paper No. 105, 10},
      ISSN = {0001-7701,1572-9532},
   MRCLASS = {53B30 (53C24 53C50)},
  MRNUMBER = {5127538},
       DOI = {10.1007/s10714-026-03609-z},
       URL = {https://doi.org/10.1007/s10714-026-03609-z},
}

@Article{galloway1989-splitting,
  author     = {Galloway, Gregory J.},
  journal    = {J. Differential Geom.},
  title      = {The {L}orentzian splitting theorem without the completeness assumption},
  year       = {1989},
  issn       = {0022-040X,1945-743X},
  number     = {2},
  pages      = {373--387},
  volume     = {29},
  fjournal   = {Journal of Differential Geometry},
  mrclass    = {53C50},
  mrnumber   = {982181},
  mrreviewer = {Steen\ Markvorsen},
  url        = {http://projecteuclid.org/euclid.jdg/1214442881},
}

@Article{ohta2014,
  author     = {Ohta, Shinichi},
  journal    = {Anal. Geom. Metr. Spaces},
  title      = {On the curvature and heat flow on {H}amiltonian systems},
  year       = {2014},
  issn       = {2299-3274},
  number     = {1},
  pages      = {81--114},
  volume     = {2},
  doi        = {10.2478/agms-2014-0003},
  fjournal   = {Analysis and Geometry in Metric Spaces},
  mrclass    = {37N05 (49Q20 53C60 58J35)},
  mrnumber   = {3208069},
  mrreviewer = {Steven\ Rosenberg},
  url        = {https://doi.org/10.2478/agms-2014-0003},
}

@Article{chrusciel-grant2012,
  author     = {Chru\'{s}ciel, Piotr T. and Grant, James D. E.},
  journal    = {Classical Quantum Gravity},
  title      = {On {L}orentzian causality with continuous metrics},
  year       = {2012},
  issn       = {0264-9381,1361-6382},
  number     = {14},
  pages      = {145001, 32},
  volume     = {29},
  doi        = {10.1088/0264-9381/29/14/145001},
  fjournal   = {Classical and Quantum Gravity},
  mrclass    = {53C50 (83A05)},
  mrnumber   = {2949547},
  mrreviewer = {Miguel\ S\'{a}nchez},
  url        = {https://doi.org/10.1088/0264-9381/29/14/145001},
}

@Article{hartman-wintner1951,
  author     = {Hartman, Philip and Wintner, Aurel},
  journal    = {Amer. J. Math.},
  title      = {On the problems of geodesics in the small},
  year       = {1951},
  issn       = {0002-9327,1080-6377},
  pages      = {132--148},
  volume     = {73},
  doi        = {10.2307/2372166},
  fjournal   = {American Journal of Mathematics},
  mrclass    = {53.0X},
  mrnumber   = {40740},
  mrreviewer = {A.\ Fialkow},
  url        = {https://doi.org/10.2307/2372166},
}

@Article{newman1990,
  author     = {Newman, Richard P. A. C.},
  journal    = {J. Differential Geom.},
  title      = {A proof of the splitting conjecture of {S}.-{T}.\ {Y}au},
  year       = {1990},
  issn       = {0022-040X,1945-743X},
  number     = {1},
  pages      = {163--184},
  volume     = {31},
  fjournal   = {Journal of Differential Geometry},
  mrclass    = {53C50},
  mrnumber   = {1030669},
  mrreviewer = {Paul\ E.\ Ehrlich},
  url        = {http://projecteuclid.org/euclid.jdg/1214444093},
}

@Book{gilbarg-trudinger2001,
  author    = {Gilbarg, David and Trudinger, Neil S.},
  publisher = {Springer-Verlag, Berlin},
  title     = {Elliptic partial differential equations of second order},
  year      = {2001},
  isbn      = {3-540-41160-7},
  note      = {Reprint of the 1998 edition},
  series    = {Classics in Mathematics},
  mrclass   = {35-02 (35Jxx)},
  mrnumber  = {1814364},
  pages     = {xiv+517},
}

@article{braun2024+,
    AUTHOR = {Braun, Mathias},
     TITLE = {Exact d'{A}lembertian for {L}orentz distance functions},
   JOURNAL = {Calc. Var. Partial Differential Equations},
  FJOURNAL = {Calculus of Variations and Partial Differential Equations},
    VOLUME = {65},
      YEAR = {2026},
    NUMBER = {9},
     PAGES = {Paper No. 262, 87},
      ISSN = {0944-2669,1432-0835},
   MRCLASS = {53C50 (28A50 35J92 53C21)},
  MRNUMBER = {5129163},
       DOI = {10.1007/s00526-026-03429-8},
       URL = {https://doi.org/10.1007/s00526-026-03429-8},
}

@article{braun-gigli-mccann-ohanyan-samann+,
    AUTHOR = {Braun, Mathias and Gigli, Nicola and McCann, Robert J. and
              Ohanyan, Argam and S\"amann, Clemens},
     TITLE = {An elliptic proof of the splitting theorems from {L}orentzian
              geometry},
   JOURNAL = {Commun. Am. Math. Soc.},
  FJOURNAL = {Communications of the American Mathematical Society},
    VOLUME = {6},
      YEAR = {2026},
     PAGES = {1087--1112},
      ISSN = {2692-3688},
   MRCLASS = {83C75 (35J92 35Q75 49Q22 51K10 53C50 58J05)},
  MRNUMBER = {5114214},
       DOI = {10.1090/cams/77},
       URL = {https://doi.org/10.1090/cams/77},
}

@Misc{braun-gigli-mccann-ohanyan-samann++,
  author = {Braun, Mathias and Gigli, Nicola and McCann, Robert J. and Ohanyan, Argam and S\"amann, Clemens},  
  title  = {{A} {L}orentzian splitting theorem for continuously differentiable metrics and weights},
  year   = {Ann. PDE, to appear},
}

@article{campanato1963,
    AUTHOR = {Campanato, S.},
     TITLE = {Propriet\`a{} di h\"olderianit\`a{} di alcune classi di
              funzioni},
   JOURNAL = {Ann. Scuola Norm. Sup. Pisa Cl. Sci. (3)},
  FJOURNAL = {Annali della Scuola Normale Superiore di Pisa. Classe di
              Scienze. Serie III},
    VOLUME = {17},
      YEAR = {1963},
     PAGES = {175--188},
      ISSN = {0391-173X},
   MRCLASS = {46.38},
  MRNUMBER = {156188},
MRREVIEWER = {C.\ B.\ Morrey, Jr.},
}

@book{evans-gariepy1992,
    AUTHOR = {Evans, Lawrence C. and Gariepy, Ronald F.},
     TITLE = {Measure theory and fine properties of functions},
    SERIES = {Studies in Advanced Mathematics},
 PUBLISHER = {CRC Press, Boca Raton, FL},
      YEAR = {1992},
     PAGES = {viii+268},
      ISBN = {0-8493-7157-0},
   MRCLASS = {28-02 (26-02 26Bxx 46E35)},
  MRNUMBER = {1158660},
MRREVIEWER = {R.\ G.\ Bartle},
}

@article{meyers1964,
    AUTHOR = {Meyers, Norman G.},
     TITLE = {Mean oscillation over cubes and {H}\"older continuity},
   JOURNAL = {Proc. Amer. Math. Soc.},
  FJOURNAL = {Proceedings of the American Mathematical Society},
    VOLUME = {15},
      YEAR = {1964},
     PAGES = {717--721},
      ISSN = {0002-9939,1088-6826},
   MRCLASS = {26.55},
  MRNUMBER = {168712},
MRREVIEWER = {L.\ C.\ Young},
       DOI = {10.2307/2034586},
       URL = {https://doi.org/10.2307/2034586},
}

@Misc{mccann+,
  author = {McCann, Robert J.},
  title  = {{T}rading linearity for ellipticity: a nonsmooth approach to {E}instein's theory of gravity and the {L}orentzian splitting theorems},
  year   = {Proceedings of the Forward From the Fields Medal 2024, to appear},
}

@Article{choquet-bruhat1952,
  author     = {Choquet-Bruhat, Y.},
  journal    = {Acta Math.},
  title      = {Th\'eor\`eme d'existence pour certains syst\`emes d'\'equations aux d\'eriv\'ees partielles non lin\'eaires},
  year       = {1952},
  issn       = {0001-5962,1871-2509},
  pages      = {141--225},
  volume     = {88},
  doi        = {10.1007/BF02392131},
  fjournal   = {Acta Mathematica},
  mrclass    = {36.0X},
  mrnumber   = {53338},
  mrreviewer = {F.\ John},
  url        = {https://doi.org/10.1007/BF02392131},
}

@Article{hawking-king-mccarthy1976,
  author     = {Hawking, S. W. and King, A. R. and McCarthy, P. J.},
  journal    = {J. Mathematical Phys.},
  title      = {A new topology for curved space-time which incorporates the causal, differential, and conformal structures},
  year       = {1976},
  issn       = {0022-2488,1089-7658},
  number     = {2},
  pages      = {174--181},
  volume     = {17},
  doi        = {10.1063/1.522874},
  fjournal   = {Journal of Mathematical Physics},
  mrclass    = {83.57},
  mrnumber   = {395736},
  mrreviewer = {S.\ Nanda},
  url        = {https://doi.org/10.1063/1.522874},
}

@article{braun2025,
    AUTHOR = {Braun, Mathias},
     TITLE = {New perspectives on the d'{A}lembertian from general
              relativity. {A}n invitation},
   JOURNAL = {Indag. Math. (N.S.)},
  FJOURNAL = {Koninklijke Nederlandse Akademie van Wetenschappen.
              Indagationes Mathematicae. New Series},
    VOLUME = {37},
      YEAR = {2026},
    NUMBER = {3},
     PAGES = {843--891},
      ISSN = {0019-3577,1872-6100},
   MRCLASS = {28A50 (35J60 49Q22 52A20 53C80 83C75)},
  MRNUMBER = {5063898},
       DOI = {10.1016/j.indag.2025.05.002},
       URL = {https://doi.org/10.1016/j.indag.2025.05.002},
}

@book{bauschke-combettes2017,
    AUTHOR = {Bauschke, Heinz H. and Combettes, Patrick L.},
     TITLE = {Convex analysis and monotone operator theory in {H}ilbert
              spaces},
    SERIES = {CMS Books in Mathematics/Ouvrages de Math\'ematiques de la
              SMC},
   EDITION = {Second},
      NOTE = {With a foreword by H\'edy Attouch},
 PUBLISHER = {Springer, Cham},
      YEAR = {2017},
     PAGES = {xix+619},
      ISBN = {978-3-319-48310-8; 978-3-319-48311-5},
   MRCLASS = {49-02 (41A65 46B20 46C05 47H05 90C25)},
  MRNUMBER = {3616647},
       DOI = {10.1007/978-3-319-48311-5},
       URL = {https://doi.org/10.1007/978-3-319-48311-5},
}

@article{minty1962,
    AUTHOR = {Minty, George J.},
     TITLE = {Monotone (nonlinear) operators in {H}ilbert space},
   JOURNAL = {Duke Math. J.},
  FJOURNAL = {Duke Mathematical Journal},
    VOLUME = {29},
      YEAR = {1962},
     PAGES = {341--346},
      ISSN = {0012-7094,1547-7398},
   MRCLASS = {47.80},
  MRNUMBER = {169064},
MRREVIEWER = {Isaac\ V.\ Shragin},
       URL = {http://projecteuclid.org/euclid.dmj/1077470254},
}

@article{browder1963,
    AUTHOR = {Browder, Felix E.},
     TITLE = {Nonlinear elliptic boundary value problems},
   JOURNAL = {Bull. Amer. Math. Soc.},
  FJOURNAL = {Bulletin of the American Mathematical Society},
    VOLUME = {69},
      YEAR = {1963},
     PAGES = {862--874},
      ISSN = {0002-9904},
   MRCLASS = {35.45},
  MRNUMBER = {156116},
MRREVIEWER = {J.\ Peetre},
       DOI = {10.1090/S0002-9904-1963-11068-X},
       URL = {https://doi.org/10.1090/S0002-9904-1963-11068-X},
}

@book{zeidler1990b,
    AUTHOR = {Zeidler, Eberhard},
     TITLE = {Nonlinear functional analysis and its applications. {II}/{B}},
      NOTE = {Nonlinear monotone operators,
              Translated from the German by the author and Leo F. Boron},
 PUBLISHER = {Springer-Verlag, New York},
      YEAR = {1990},
     PAGES = {i--xvi and 469--1202},
      ISBN = {0-387-97167-X},
   MRCLASS = {47-02 (35-01 35J60 47Hxx 58-01 65Jxx)},
  MRNUMBER = {1033498},
MRREVIEWER = {Jean\ Mawhin},
       DOI = {10.1007/978-1-4612-0985-0},
       URL = {https://doi.org/10.1007/978-1-4612-0985-0},
}

@Misc{ryborz2026+,
  author = {Ryborz, Vanessa},
  title  = {{I}nfinitesimal {M}inkowskianity for manifolds with continuous {L}orentzian metrics},
  year   = {Preprint, arXiv:2604.22420},
}

@book{giaquinta-martinazzi2012,
    AUTHOR = {Giaquinta, Mariano and Martinazzi, Luca},
     TITLE = {An introduction to the regularity theory for elliptic systems,
              harmonic maps and minimal graphs},
    SERIES = {Appunti. Scuola Normale Superiore di Pisa (Nuova Serie)
              [Lecture Notes. Scuola Normale Superiore di Pisa (New
              Series)]},
    VOLUME = {11},
   EDITION = {Second},
 PUBLISHER = {Edizioni della Normale, Pisa},
      YEAR = {2012},
     PAGES = {xiv+366},
      ISBN = {978-88-7642-442-7; 978-88-7642-443-4},
   MRCLASS = {35-02 (35B65 35J20 35J60 58E20)},
  MRNUMBER = {3099262},
       DOI = {10.1007/978-88-7642-443-4},
       URL = {https://doi.org/10.1007/978-88-7642-443-4},
}

@Misc{ohanyan2026+,
  author = {Ohanyan, Argam},
  title  = {Cosmological rigidity via $p$-harmonic temporal functions},
  year   = {Preprint, arXiv:2611.01456},
}

@Article{myers-steenrod1939,
  author   = {Myers, S. B. and Steenrod, N. E.},
  journal  = {Ann. of Math. (2)},
  title    = {The group of isometries of a {R}iemannian manifold},
  year     = {1939},
  issn     = {0003-486X,1939-8980},
  number   = {2},
  pages    = {400--416},
  volume   = {40},
  doi      = {10.2307/1968928},
  fjournal = {Annals of Mathematics. Second Series},
  mrclass  = {99-04},
  mrnumber = {1503467},
  url      = {https://doi.org/10.2307/1968928},
}

@Article{calabi-hartman1970,
  author     = {Calabi, Eugenio and Hartman, Philip},
  journal    = {Duke Math. J.},
  title      = {On the smoothness of isometries},
  year       = {1970},
  issn       = {0012-7094,1547-7398},
  pages      = {741--750},
  volume     = {37},
  fjournal   = {Duke Mathematical Journal},
  mrclass    = {53.74},
  mrnumber   = {283727},
  mrreviewer = {N.\ J.\ Hicks},
  url        = {http://projecteuclid.org/euclid.dmj/1077379303},
}

@Article{palais1957,
  author     = {Palais, Richard S.},
  journal    = {Proc. Amer. Math. Soc.},
  title      = {On the differentiability of isometries},
  year       = {1957},
  issn       = {0002-9939,1088-6826},
  pages      = {805--807},
  volume     = {8},
  doi        = {10.2307/2033302},
  fjournal   = {Proceedings of the American Mathematical Society},
  mrclass    = {53.2X},
  mrnumber   = {88000},
  mrreviewer = {K.\ Krickeberg},
  url        = {https://doi.org/10.2307/2033302},
}

\end{document}